\documentclass[cupthm]{CUP-JNL-JMJ}

\usepackage{latexsym}
\usepackage{graphicx}
\usepackage{longtable}
\usepackage{booktabs}
\usepackage{multicol,multirow}
\usepackage{amsmath,amssymb,amsfonts}
\usepackage{mathrsfs}
\usepackage{amsthm}
\usepackage{rotating}
\usepackage{appendix}
\usepackage[numbers]{natbib}
\usepackage{tikz}
\usetikzlibrary{cd}

\numberwithin{equation}{section}

\articletype{RESEARCH ARTICLE}
\jname{J. Inst. Math. Jussieu}
\jyear{2020}
\jvol{1}
\jdoi{10.1017/S1474748020000237}

\theoremstyle{cupplain}
\newtheorem{theorem}{Theorem}[section]
\newtheorem{lemma}[theorem]{Lemma}
\newtheorem{corollary}[theorem]{Corollary}

\theoremstyle{cupdefinition}

\theoremstyle{cupremark}

\theoremstyle{cupproof}
\newtheorem{proof}{Proof}

\begin{document}

\begin{Frontmatter}

\title[The geometric expansion of the coarse trace formula of $\rm{GL}(4)$]{The geometric expansion of the coarse trace formula of $\rm{GL}(4)$\thanks{Research supported by National natural Science
Foundation of PR China (No. 12071326)}}

\author{Haoyang Wang, Xinghua Cui, Weilong Zhang, Zhifeng Peng*\thanks{*Corresponding author.}}

\address{\orgdiv{School of Mathematical Science.} \orgname{Soochow University}, \orgaddress{\state{JiangSu}, \city{Suzhou, P.R.China}} \email{20244207002@stu.suda.edu.cn, 20234007004@stu.suda.edu.cn, 20254207011@stu.suda.edu.cn, zfpeng@suda.edu.cn}}

\received{12 June 2026}

\maketitle

\abstract{We study the non-semisimple terms in the geometric side of Arthur's trace formula of $\rm{GL}(4)$ for the further study of the invariant trace formula. After modifying the zeta function and Arthur's definition of the weighted orbital integrals, we compare the two ways to calculate the distributions $J^T(f)$ of the geometric side, we can obtain that the distributions can be written as the linear combanition of the weighted orbital integrals. In particular, we get the coefficients $a^M(S,\gamma)$, for all $\gamma\in \rm{GL}(4,F)$, expressed by the measure and zeta functions.}


\keywords[\textup{2020} Mathematics subject classification]{11F72}


\end{Frontmatter}

\tableofcontents
\section{Introduction}
James Arthur, in a monumental series of papers spanning several decades, developed the Arthur trace formula for general connected reductive groups $G$. Arthur's insight was to introduce a subtle, non-trivial truncation operator $\Lambda^T$ that acts on the kernel function of the right regular representation on $L^2(G(F)\backslash G(\mathbb{A})^1)$, for the number field $F$ and the Adele ring $\mathbb{A}$ of $F$. By carefully analyzing the asymptotic behavior of Eisenstein series along the roots of the parabolic subgroups, Arthur proved that the truncated kernel is integrable, thereby yielding a canonical identity:
\begin{eqnarray}
    J_{\rm{geom}}(f)=J_{\rm{spec}}(f).
\end{eqnarray}
Here, both sides are distributions on the test function space $C_c^\infty(G(\mathbb{A})^1)$. The geometric side $J_{\rm{geom}}(f)$ is a sum of terms indexed by equivalance classes of elements in $G(F)$, which can be further decomposed into weighted orbital integrals. The spectral side $J_{\rm{spec}}(f)$ contains not only the discrete spectrum (which encodes the arithmetic data of interest, such as Galois representations) but also contributions from the continuous spectrum expressed via weighted traces and intertwining operators.
In subsequent developments, Arthur refined this identity into the invariant trace formula (in \cite{A4},\cite{A7},\cite{A8}), ensuring that each term on both sides is invariant under conjugation. This refinement paved the way for the stabilization of the trace formula that expresses the geometric side in terms of stable conjugacy classes on $G$ and its endoscopic datum.

In coarse trace formula, the distributions of unramified orbits $\mathfrak{o}$ is \( J_{\mathfrak{o}}^T(f) \), which is 
\[
\rm{vol}\bigl( A_0(\mathbb{R})^0 G_{\mu_{\mathfrak{o}}}(F) \setminus G_{\mu_{\mathfrak{o}}}(\mathbb{A}) \bigr)
\int_{G_{\mu_{\mathfrak{o}}}(\mathbb{A}) \setminus G(\mathbb{A})} f(x^{-1} \mu x) \, v(x, T) \, dx,
\]
where  
\[
v(x, T) = \int_{Z(\mathbb{R})^0 \setminus A_\mathfrak{o}(\mathbb{R})^0}
\left\{
\sum_{P} (-1)^{\dim(A/Z)} \sum_{s \in \Omega (\mathfrak{a}_\mathfrak{o}, P)}
\hat{\tau}_P\bigl( H_0(w_s a x) - T \bigr)
\right\} da,
\] and $\mu_\mathfrak{o}$ depends on $\mathfrak{o}$. Sometimes $v(x,T)$ is written as $v_{M_\mathfrak{o}}(x)$ or $v_{M}(x)$ if we take $T=T_0$ ($T_0$ is dependent on the group $G$, and the existence is proved in \cite{A4}).

For ramfied orbits $\mathfrak{o}$, we could only write $J_{\mathfrak{o}}^T(f)$ as the integral over $G(F)\backslash G(\mathbb{A})^1$ of
\[\sum_P (-1)^{\dim(A/Z)} \cdot \sum_{\delta \in P(F) \setminus G(F)} \sum_{\gamma\in M(F)\cap\mathfrak{o}}\int_{N(
\mathbb{A})}f(x^{-1}\delta^{-1}\gamma n\delta x)dn\hat{\tau}_P(H_0(\delta x) - T).\]

To obtain the invariant trace formula, we need to perform additional operations on both sides of the trace formula. Arthur gives a way to deal with the convex hull $v_M(x)$ in order to get the invariance, but one question is from the distributions of ramified terms, they do not have convex hull. Arthur in \cite{A6} proved that the distributions of ramified terms can be written as the limit of the the sum of weighted orbital integrals:
\begin{eqnarray}
    J^T(f)=\sum_{M\in\mathcal{L}}|W_0^M||W_0^G|^{-1}\sum_{\gamma\in(M(F))_{M,S}}a^M(S,\gamma)J_M^T(\gamma,f),
\end{eqnarray}
where $J^T(f)$ equals the integral over $G(F)\backslash G(\mathbb{A})^1$ of 
\begin{eqnarray}\label{1.1}
\\\notag
    \sum_{P\supset P_0}(-1)^{\rm{dim}\, Z\backslash A_P}\sum_{\delta\in P(F)\backslash G(F)}\sum_{\gamma\in M(F)\cap\mathfrak{o}}\int_{N(
\mathbb{A})}f(x^{-1}\delta^{-1}\gamma n\delta x)dn\hat{\tau}_P(H_P(\delta x)-T_P),
\end{eqnarray}
and
\begin{eqnarray*}
    J^T_M(\gamma, f) = |D(\gamma)|^{1/2} \int_{G_\gamma(\mathbb{A}) \setminus G(\mathbb{A})} f(x^{-1}\gamma x) \, \nu_M(x) \, dx.
\end{eqnarray*}
Therefore, we can directly work with the convex hull of each term on the geometric side of the trace formula to obtain the invariance.

For any element $\gamma\in \rm{GL}(4,F)$, write $\sigma$ as the semisimple part of the Jordan decomposition of $\gamma$, define \begin{eqnarray*}
\epsilon^M(\sigma)=\begin{cases}
        1,\quad &if\  A_{M_{\sigma}}=A_M,\\
        0,\quad &else,
    \end{cases}
\end{eqnarray*}for any Levi subgroup $M$.

Arthur has proved the descend formula:
\begin{eqnarray*}
a^M(S,\gamma)=\epsilon^M(\sigma)|\iota^M(\sigma)|^{-1}\sum_{\{u|\sigma u\sim \gamma\}}a^{M_{\sigma}}(S,u),
\end{eqnarray*}
where the sum is taken over all $u\in (\mathcal{U}_{M_\sigma}(F))_{M_{\sigma,S}}$ such that $\sigma u$ is $(M,S)$-equivalent to $\gamma$. $\iota^M(\sigma)$ is  $|M_{\sigma}^0(F)\backslash M_{\sigma}(F)|$, that is, the number of the connected components. In the case of $\rm{GL}(n)$, it always equals $1$.

It means that the other elements can not conjugated to a unipotent elements by Jordan decomposition, we can descent to the unipotent case, though some $a^M(S,\gamma)$ may equals $0$. Thus we just need to calculate the coefficient $a^M(S,\gamma)$ for $\gamma$ is the unipotent case.

The coefficients $a^M(S,\gamma)$ have determined only in special cases. For semisimple elements $\gamma$, Arthur gives $a^M(S,\gamma)$ in \cite{A6}, for $M$ of $F$-rank one, Hoffmann gives all $a^M(S,\gamma)$ in \cite{Ho} and for $M=\rm{GL}(3)$, Flicker and Matz give all $a^M(S,\gamma)$ in \cite{F} and \cite{M}, for $G=\rm{GSp}(2)$ and $\rm{Sp}(2)$, Hoffmann and Wakatsuki give all $a^M(S,\gamma)$ in \cite{HW}. Then for two different cases of unipotent conjugacy classes(that is, the semisimple part of the Jordan decomposition is trivial), Chaudouard deals with the weighted orbital integrals of $\rm{GL}(n)$ in \cite{C1} and \cite{C2}, but he did not give all the coefficients $a^M(S,\gamma)$ here.

These previous works compare the two methods for calculating the geometric side of the trace formula to obtain the coefficient $a^M(S,\gamma)$, as we will do in this paper.

Write $u(n_1,n_2,n_3,n_4,n_5,n_6)$ as the matrix
\begin{eqnarray*}
    \left(\begin{matrix}
        1&n_1&n_2&n_3\\
        &1&n_4&n_5\\
        &&1&n_6\\
        &&&1
    \end{matrix}\right).
\end{eqnarray*}

The main result of this article is that we obtain all the coefficients $a^M(S,\gamma)$ in the case of $\rm{GL}(4)$. The approach we shall use replaces the invariance argument
of \cite{A3} and removes the restriction on $S$. For the unipotent conjugacy classes, we have proved the following theorems. For the specific notation, the reader is referred to the text.

Let $z\in Z(F)$, which is the center of $G(F)$.

For the conjugacy class of $zu(0,0,1,0,0,0)$ in $G(F)$, we have 
\begin{theorem}(Theorem \ref{theorem1})
    For $\mathfrak{o}=\mathfrak{o}_{31,z}$, $z\in Z(F)$ and $f\in C_c^\infty(G(\mathbb{A})^1)$, $J_{\mathfrak{o}_{31,z}}(f)$ equals
    \[\rm{vol}_{M_{31}}J_{M_{31}}(1,f)+\rm{vol}_{M_{31}}\cdot \frac{\frac{d}{ds}\zeta_F^S(s)|_{s=3}}{\zeta_F^S(3)}J_G(u(0,0,1,0,0,0),f),\]that is \[a^{M_{31}}(S,1)=\rm{vol}_{M_{31}},\quad a^{G}(S,u(0,0,1,0,0,0))=\rm{vol}_{M_{31}}\cdot \frac{\frac{d}{ds}\zeta_F^S(s)|_{s=3}}{\zeta_F^S(3)}.\]
\end{theorem}
For the conjugacy class of $u(1,0,0,0,0,1)$ in $G(F)$, we have 
\begin{theorem}(Theorem \ref{theorem2})
    For $\mathfrak{o}=\mathfrak{o}_{22,z}$, $z\in Z(F)$ and $f\in C_c^\infty(G(\mathbb{A})^1)$, $J_{\mathfrak{o}_{22,z}}(f)$ equals
    \[\rm{vol}_{M_{22}}J_{M_{22}}(1,f)+\rm{vol}_{M_{22}}\cdot \mathfrak{c}_{2,F}(S)J_G(u(0,1,0,0,1,0),f),\]that is \[a^{M_{22}}(S,1)=\rm{vol}_{M_{22}},\quad a^{G}(S,u(0,1,0,0,1,0))=\rm{vol}_{M_{22}}\cdot \mathfrak{c}_{2,F}(S).\]
\end{theorem}
For the conjugacy class of $u(0,1,0,0,0,1)$ in $G(F)$, we have 
\begin{theorem}(Theorem \ref{theorem3})
    For $\mathfrak{o}=\mathfrak{o}_{211,z}$, $z\in Z(F)$ and $f\in C_c^\infty(G(\mathbb{A})^1)$, $J_{\mathfrak{o}_{211,z}}(f)$ equals
    \[\rm{vol}_{M_{31}}J_{M_{31}}(1,f)+\rm{vol}_{M_{31}}\cdot \frac{\frac{d}{ds}\zeta_F^S(s)|_{s=3}}{\zeta_F^S(3)}J_G(u(0,0,1,0,0,0),f),\]that is 
    \begin{eqnarray*}
        &a^{M_{211}}(S,1)=\rm{vol}_{M_{211}},\quad a^{M_{22}}(S,u(0,0,0,0,0,1))=\rm{vol}_{M_{211}}\frac{\mathfrak{c}_F(S)}{2c_F},\\
        &a^{M_{31}}(S,u(0,1,0,0,0,0))=\rm{vol}_{M_{211}}\frac{\frac{d}{ds}\zeta_F^S(s)|_{s=2}}{\zeta_F^S(2)},\\
        &a^G(S,u(0,1,0,0,0,1))=\rm{vol}_{M_{211}}\frac{\frac{d^2}{ds^2}\zeta_F^S(s)|_{s=2}}{\zeta_F^S(2)}.
    \end{eqnarray*}
\end{theorem}
For the conjugacy class of $u(1,0,0,1,0,1)$ in $G(F)$, we have 
\begin{theorem}(Theorem \ref{theorem4})
    For $\mathfrak{o}=\mathfrak{o}_{1,z}$, $z\in Z(F)$ and $f\in C_c^\infty(G(\mathbb{A})^1)$, $J_{\mathfrak{o}_{1,z}}(f)$ equals
    \[\rm{vol}_{M_{0}}J_{M_{0}}(1,f)+\rm{vol}_{M_{22}}\cdot \mathfrak{c}_{2,F}(S)J_G(u(0,1,0,0,1,0),f),\]that is 
    \begin{eqnarray*}
        &a^{M_{0}}(S,1)=\rm{vol}_{M_{0}},\quad a^{M_{211}}(S,u(1,0,0,0,0,0))=\rm{vol}_{M_0}\frac{\mathfrak{c}_F(S)}{2c_F},\\
        &a^{M_{31}}(S,u(1,0,0,1,0,0))=\frac{\rm{vol}_{M_0}}{3c_F^2}(\mathfrak{c}_F(S)^2+\mathfrak{c}'_F(S)c_F^S),\\
        &a^{M_{22}}(S,u(1,0,0,0,0,1))=\frac{\rm{vol}_{M_0}}{4c_F^2}\mathfrak{c}_F^2(S),\\
        &a^{G}(S,u(0,1,0,0,1,0))=\frac{\rm{vol}_{M_0}}{4c_F^3}(\mathfrak{c}_F^3(S)+3\mathfrak{c}'_F(S)\mathfrak{c}_F(S)c_F^S+\mathfrak{c}''_F(S)(c_F^S)^2).
    \end{eqnarray*}
\end{theorem}
And for the elements whose semisimple parts are not trivial, we deal with this case according to Arthur's results in section \ref{uniddd}.

Next comes the organization of the chapters in this article.

In Section~\ref{2}, we collect the necessary preliminaries and fix notation. In Section~\ref{3}, we recall the main results of Arthur's weighted orbital integrals and try to give a modification for the further use. In Section~\ref{4} and \ref{5}, we introduce the definition of zeta integrals and give a modification to make it cater to the calculation of the geometric terms of the coarse trace formula of $\rm{GL}(4)$. In Section~\ref{6}, we describe the mean value of $\rm{GL}(3)$ which we shall use later. In Section~\ref{7}, we divide the trace formula of $\rm{GL}(4)$ according to our previous work of coarse trace of $\rm{GL}(4)$ in \cite{WCP} and give a explicit expression of that trace formula. In Section~\ref{8}, we modify the arthur's method to simplify the calculation by using the idea of cancelling the poles of Laurent expansion. Finally, we use the result of Arthur's descent formula to give all the coefficient $a^M(S,\gamma)$.

\section{Preliminaries}\label{2}
\subsection{Structure of parabolic subgroups}
For global field $F$ and connected linear algebraic $F$-group $G=\rm{GL}(4)$. Let $S$ be the finite set of the places of $F$ which contains archimedean places and $2$.  Write $\Sigma$ as the set of places.  Then $G(\mathbb{A})$ is the restricted direct product over all valuations $v$ of the groups $G(F_v)$. We fix the maximal compact subgroup
$$K= \prod_{v \in \Sigma} K_v \quad \text{where} \quad K_v = \begin{cases} \mathrm{U}(4) & \text{if } v \in \Sigma_{\mathbb{C}}, \\ \mathrm{O}(4) & \text{if } v \in \Sigma_{\mathbb{R}}, \\ \mathrm{GL}(4, \mathcal{O}_v) & \text{if } v \in \Sigma_{\mathrm{fin}}. \end{cases}$$
Fix the standard maximal torus $T$ consisting of diagonal matrices, and the standard minimal parabolic subgroup (or Borel subgroup) $B=P_
0$ consisting of upper triangular matrices in $G$.

Denote $X(G)_F$ by the group of $F$-rational characters of $G$, $A_G$ by the $F$-split part of the center of $G$. Define the vector space $\mathfrak{a}_G=\rm{Hom}_\mathbb{Z}(X(G)_F,\mathbb{R})$. We can similarly define these notation by  replacing $G$ by any standard Levi subgroup $M$.

Then if $m=\Pi_vm_v$ lies in $M(\mathbb{A})$, we define a vector $H_M(m)=H_P(m)$ in $\mathfrak{a}_P$ by 
\[\rm{exp}(\langle H_M(m),\chi\rangle)=|\chi(m)|=\Pi_v|\chi(m_v)|_v.\]
Write $M_P(\mathbb{A})^1$ be the kernel of the homomorphism $H_M$, then \[M_P(\mathbb{A})=M_P(\mathbb{A})^1A_{M,\infty}^+,\] where $A_{M,\infty}$ is the infinity place of $A_M$, $A_{M,\infty}^+$ is the connected component of $1$ in $A_{M,\infty}$. By the Langlands decomposition, $G(\mathbb{A})=N_P(\mathbb{A})M_P(\mathbb{A})K$, we can then write any $x\in G(\mathbb{A})$ as $nmk$, $n\in N_P(\mathbb{A})$, $m\in M_P(\mathbb{A})$, $k\in K$.

Set $\Omega$ be the restricted Weyl group of $(G,A_0)$. $\Omega$ acts on the dual space of $\mathfrak{a}_0$. The roots of $(P,A)$, which are the elements in $X(A)_F$ obtained by decomposing the Lie algebra of $N_P$.

We also set
$$P_{31} = \left\{ \begin{pmatrix} * & * & * & * \\ * & * & * & * \\ * & * & * & * \\ 0 & 0 & 0 & * \end{pmatrix} \right\}, \quad P_{22} = \left\{ \begin{pmatrix} * & * & * & * \\ * & * & * & * \\ 0 & 0 & * & * \\ 0 & 0 & * & * \end{pmatrix} \right\}, \quad P_{13} = \left\{ \begin{pmatrix} * & * & * & * \\ 0 & * & * & * \\ 0 & * & * & * \\ 0 & * & * & * \end{pmatrix} \right\} ,$$
\[P_{211} = \left\{ \begin{pmatrix} * & * & * & * \\ * & * & * & * \\ 0 & 0 & * & * \\ 0 & 0 & 0 & * \end{pmatrix} \right\},\quad P_{121} = \left\{ \begin{pmatrix} * & * & * & * \\ 0 & * & * & * \\ 0 & * & * & * \\ 0 & 0 & 0 & * \end{pmatrix} \right\},\quad P_{112} = \left\{ \begin{pmatrix} * & * & * & * \\ 0 & * & * & * \\ 0 & 0 & * & * \\ 0 & 0 & * & * \end{pmatrix}\right\}.\]

Then, the standard parabolic subgroups are 
$$\mathcal{F}=\{ P \in \mathcal{P} \mid P \supset P_0 \} = \{ P_0, \, P_{211}, \, P_{121}, \, P_{112}, \, P_{31}, \, P_{22}, \, P_{13}, \, G \}.$$
The group $\Omega$ acts on $\mathbb{R}^4$ as
\begin{equation*}
    s(x_1, x_2, x_3, x_4) = (x_{s(1)}, x_{s(2)}, x_{s(3)}, x_{s(4)}) \quad \text{for} \quad s \in W_0.
\end{equation*}
Let $(ij)$ denote the transposition of $i$ and $j$. We define the homomorphism $W_0 \rightarrow \mathrm{GL}(4)$ by
\begin{equation*}
    (12) \mapsto \begin{pmatrix} 0 & 1 & 0 & 0 \\ 1 & 0 & 0 & 0 \\ 0 & 0 & 1 & 0 \\ 0 & 0 & 0 & 1 \end{pmatrix}, \quad 
    (23) \mapsto \begin{pmatrix} 1 & 0 & 0 & 0 \\ 0 & 0 & 1 & 0 \\ 0 & 1 & 0 & 0 \\ 0 & 0 & 0 & 1 \end{pmatrix}, \quad 
    (34) \mapsto \begin{pmatrix} 1 & 0 & 0 & 0 \\ 0 & 1 & 0 & 0 \\ 0 & 0 & 0 & 1 \\ 0 & 0 & 1 & 0 \end{pmatrix}.
\end{equation*}
\subsection{Lie algebras and roots}
Suppose that $P_1$ and $P_2$ are two standard parabolic subgroups, with $P_1 \subset P_2$. There are then $\mathbb{Q}$-rational embeddings
\begin{equation*}
    A_{P_2} \subset A_{P_1} \subset M_{P_1} \subset M_{P_2}.
\end{equation*}
The restriction homomorphism $X(M_{P_2})_\mathbb{Q} \rightarrow X(M_{P_1})_\mathbb{Q}$ is injective. It provides a linear injection $\mathfrak{a}_{P_2}^* \hookrightarrow \mathfrak{a}_{P_1}^*$ and a dual linear surjection $\mathfrak{a}_{P_1} \twoheadrightarrow \mathfrak{a}_{P_2}$. We write $\mathfrak{a}_{P_1}^{P_2} \subset \mathfrak{a}_{P_1}$ for the kernel of the latter mapping. The restriction homomorphism $X(A_{P_1})_\mathbb{Q} \rightarrow X(A_{P_2})_\mathbb{Q}$ is surjective, and extends to a surjective mapping from $X(A_{P_1})_\mathbb{Q} \otimes \mathbb{R}$ to $X(A_{P_2})_\mathbb{Q} \otimes \mathbb{R}$. It thus provides a linear surjection $\mathfrak{a}_{P_1}^* \twoheadrightarrow \mathfrak{a}_{P_2}^*$, and a dual linear injection $\mathfrak{a}_{P_2} \hookrightarrow \mathfrak{a}_{P_1}$. The four linear mappings yield split exact sequences
\begin{equation*}
\begin{tikzcd}
    0 \arrow[r] & \mathfrak{a}_{P_2}^* \arrow[r, shift left=0.5ex] & \mathfrak{a}_{P_1}^* \arrow[l, shift left=0.5ex] \arrow[r] & \mathfrak{a}_{P_1}^* / \mathfrak{a}_{P_2}^* \arrow[r] & 0
\end{tikzcd}
\end{equation*}
and
\begin{equation*}
\begin{tikzcd}
    0 \arrow[r] & \mathfrak{a}_{P_1}^{P_2} \arrow[r] & \mathfrak{a}_{P_1} \arrow[r, shift left=0.5ex] & \mathfrak{a}_{P_2} \arrow[l, shift left=0.5ex] \arrow[r] & 0
\end{tikzcd}
\end{equation*}
of real vector spaces. We may therefore write
\begin{equation*}
    \mathfrak{a}_{P_1} = \mathfrak{a}_{P_2} \oplus \mathfrak{a}_{P_1}^{P_2}
\end{equation*}
and
\begin{equation*}
    \mathfrak{a}_{P_1}^* = \mathfrak{a}_{P_2}^* \oplus (\mathfrak{a}_{P_1}^{P_2})^*.
\end{equation*}
Let $\mathfrak{g}$ denote the Lie algebra of $G$ and let $\mathrm{Ad}$ denote the adjoint representation of $G$ on $\mathfrak{g}$. For $P \in \mathcal{F}$, write $\mathfrak{a}_P^* = \mathfrak{a}_{M_P}^*$, $\mathfrak{a}_{P_1}^{P_2} = \mathfrak{a}_{M_{P_1}}^{M_{P_2}}$, and $(\mathfrak{a}_{P_1}^{P_2})^* = (\mathfrak{a}_{M_{P_1}}^{M_{P_2}})^*$. Let $\mathfrak{n}_P$ denote the Lie algebra of $N_P$. For each $\alpha \in X(A_P)_F$, a subspace $\mathfrak{n}_\alpha$ of $\mathfrak{n}_P$ is defined by
\begin{equation*}
    \mathfrak{n}_\alpha = \{X \in \mathfrak{n}_P \mid \mathrm{Ad}(a)X = \alpha(a)X, \; \forall a \in A_P\}.
\end{equation*}
Then, a finite set $\Phi_P$ is defined by
\begin{equation*}
    \Phi_P = \{\alpha \in X(A_P)_\mathbb{Q} \mid \alpha \neq 0, \; \mathfrak{n}_\alpha \neq \{0\}\}.
\end{equation*}
    It is straightforward to verify that $\mathfrak{n}_P = \bigoplus_{\alpha \in \Phi_P} \mathfrak{n}_\alpha$. Denote by $\Delta_P$ the subset of primitive roots, i.e., those $\alpha \in \Phi_P$ which are not of the form $\beta + \gamma$ with $\beta, \, \gamma \in \Phi_P$. This is a basis of $(\mathfrak{a}_P^G)^*$.

Write $\mathfrak{a}_0 = \mathfrak{a}_{M_0}$, $\mathfrak{a}_0^* = \mathfrak{a}_{M_0}^*$, $\mathfrak{a}_0^P = \mathfrak{a}_{M_0}^{M_P}$, $\mathfrak{a}_P^G = \mathfrak{a}_{M_P}^{M_G}$, $(\mathfrak{a}_0^P)^* = (\mathfrak{a}_{M_0}^{M_P})^* = (\mathfrak{a}_{M_0}^{M_P})^*$, and $\Phi_0 = \Phi_{P_0}$. Now, $\Phi_0 \cup (-\Phi_0)$ is a root system in $(\mathfrak{a}_0^G)^*$ and $\Phi_0$ is a system of positive roots. We set $\Delta_0 = \Delta_{P_0}$, i.e., $\Delta_0$ is the set of simple roots attached to $\Phi_0$. Let $\widehat{\Delta}_0$ denote the set of simple weights corresponding to $\Delta_0$.
\subsection{Measures}
Fix a Haar measure $dg$ on $G(\mathbb{A})$ and a Haar measure $da$ on $\mathfrak{a}_P$. We denote by $dn$ the Haar measure on $N_P(\mathbb{A})$ normalized by $\int_{N_P(F) \backslash N_P(\mathbb{A})} dn = 1$. Let $dk$ denote the Haar measure on $K$ normalized by $\int_{K} dk = 1$, and let $da$ denote the Haar measure on $A_{M_P,\infty}^+$ induced from $dH$. Then, there exists a unique Haar measure $d^1m$ on $M_P(\mathbb{A})^1$ such that
\begin{equation*}
    \int_{G(\mathbb{A})} f(g) \, dg = \int_{N_P(\mathbb{A})} \int_{M_P(\mathbb{A})^1} \int_{A_{M_{P,\infty}}^+} \int_{K} f(nmak) e^{-2\rho_P(H_P(a))} \, d k \, d^1 m \, d a \,d n
\end{equation*}
for any $f \in C_c^\infty(G(\mathbb{A}))$. A Haar measure $dm$ on $M_P(\mathbb{A})$ is determined by $dm = d^1m \, da$. 

For each place $v $, let $dg_v$ denote a Haar measure on $G(F_v)$, and let $dn_v$, $dm_v$, and $dk_v$ denote Haar measures on $N_P(F_v)$, $M_P(F_v)$, and $K_v$ respectively. We may assume that $dg = \prod_{v \in \Sigma} dg_v$, $dk = \prod_{v \in \Sigma} dk_v$, $dn = \prod_{v \in \Sigma} dn_v$, and $dm = \prod_{v \in \Sigma} dm_v$. Furthermore, we may also assume that
\begin{equation*}
    \int_{K_v} dk_v = 1 \quad (v \in \Sigma), \qquad \int_{N_P(F_v) \cap \mathbf{K}_v} dn_v = 1 \quad (v \in \Sigma_{\mathrm{fin}})
\end{equation*}
and
\begin{equation*}
    \int_{G(F_v)} f_v(g_v) \, d g_v = \int_{N_P(F_v)} \int_{M_P(F_v)} \int_{K_v} f_v(n_v m_v k_v) e^{-2\rho_P(H_P(m_v))} \, dk_v \, dm_v \, dn_v
\end{equation*}
for any $v \in \Sigma$ and any $f_v \in C_c^\infty(G(F_v))$. We also have the Haar measures $dn_S = \prod_{v \in S} dn_v$, $dm_S = \prod_{v \in S} dm_v$, and $dk_S = \prod_{v \in S} dk_v$ on $N(F_S)$, $M(F_S)$, and $K_S = \prod_{v \in S} K_v$ respectively.

Define a vector \[\rho_{P}=\frac{1}{2}\sum_{\alpha\in\Phi_P}(\rm{dim}\ \mathfrak{n}_\alpha)\alpha.\]  
Then the left and right Haar measure on the group $P(\mathbb{A})$ are related by \[d_lp=\delta_P(p)d_rp,\quad p\in P(\mathbb{A}),\]where $\delta_P(p)=e^{2\rho_P(H_P(p))}$ is the modular function. In particular, the group $P(\mathbb{A})$ is not unimodular if $P\neq G$.

Note that according to our work in another paper \cite{WCP}, we omit the discussion of the absolute convergence of the integrals.

\subsection{Coordinates on $\mathfrak{a}_0$}
The rational homomorphism $\chi_j : M_0 \rightarrow \mathrm{GL}(1)$ is defined by
\begin{equation*}
    \chi_j \left( \begin{pmatrix} a_1 & 0 & 0 & 0 \\ 0 & a_2 & 0 & 0 \\ 0 & 0 & a_3 & 0 \\ 0 & 0 & 0 & a_4 \end{pmatrix} \right) = a_j, \quad 1 \le j \le 4.
\end{equation*}
Then, the free abelian group $X(M_0)_F$ is \[X(M_0)_F = \{ \chi_1^{k_1} \chi_2^{k_2} \chi_3^{k_3} \chi_4^{k_4} \mid k_1, k_2, k_3, k_4 \in \mathbb{Z} \}.\] We define $H_j \in \mathfrak{a}_0$ by
\begin{equation*}
    H_j(\chi_1^{k_1} \chi_2^{k_2} \chi_3^{k_3} \chi_4^{k_4}) = k_j.
\end{equation*}
Then, we see that
\begin{equation*}
    \mathfrak{a}_0^* = \mathbb{R}\chi_1 \oplus \mathbb{R}\chi_2 \oplus \mathbb{R}\chi_3 \oplus \mathbb{R}\chi_4 \quad \text{and} \quad \mathfrak{a}_0 = \mathbb{R}H_1 \oplus \mathbb{R}H_2 \oplus \mathbb{R}H_3 \oplus \mathbb{R}H_4.
\end{equation*}

We define the inner product $\langle \, , \rangle$ on $\mathbb{R}^4$ by
$$\langle x, y \rangle = x_1 y_1 + x_2 y_2 + x_3 y_3 + x_4 y_4,$$
where $x = (x_1, x_2, x_3, x_4), \, y = (y_1, y_2, y_3, y_4) \in \mathbb{R}^4$, it is a positive definite $\Omega$-invariant bilinear form. Then, we can identify $\mathfrak{a}_0$ with $\mathfrak{a}_0^*$ by $\langle \, , \rangle$, since $x(y) = \langle x, y \rangle$ for $x \in \mathfrak{a}_0$ and $y \in \mathfrak{a}_0^*$. We set
$$e_1 = (1,0,0,0), \quad e_2 = (0,1,0,0), \quad e_3 = (0,0,1,0), \quad e_4 = (0,0,0,1).$$

The surjective map $H_{M_0} : M_0(\mathbb{A}) \to \mathfrak{a}_0$ is given by
$$H_{M_0}(m) = (\rm{log}|m_1|, \rm{log}|m_2|, \rm{log}|m_3|, \rm{log}|m_4|) \quad \text{for } m = \mathrm{diag}(m_1, m_2, m_3, m_4) \in M_0(\mathbb{A}).$$

It is straightforward to verify that
$$\mathfrak{a}_G = \mathfrak{a}_G^* = \mathbb{R}(e_1 + e_2 + e_3 + e_4),$$
$$\mathfrak{a}_0^G = (\mathfrak{a}_0^G)^* = \{(x_1, x_2, x_3, x_4) \in \mathbb{R}^4 \mid x_1 + x_2 + x_3 + x_4 = 0\},$$
$$H_G(g) = \frac{\rm{log}|\det(g)|}{4}(e_1 + e_2 + e_3 + e_4) \in \mathfrak{a}_G \quad (g \in G(\mathbb{A})).$$

Since $P_0$, $\Phi_0 = \{ e_i - e_j \mid 1 \le i < j \le 4 \}$, we have
$$\Delta_0 = \{ \alpha_1=e_1 - e_2, \, \alpha_2=e_2 - e_3, \, \alpha_3=e_3 - e_4 \} \quad$$ and $$\hat{\Delta}_0 = \{ \hat{\alpha}_1, \, \hat{\alpha}_2, \, \hat{\alpha}_3\}$$
where we set
$$\hat{\alpha}_1 = \frac{3}{4}e_1 - \frac{1}{4}(e_2 + e_3 + e_4),$$
$$\hat{\alpha}_2= \frac{1}{2}(e_1 + e_2) - \frac{1}{2}(e_3 + e_4),$$
$$\hat{\alpha}_3= \frac{1}{4}(e_1 + e_2 + e_3) - \frac{3}{4}e_4.$$
We then can get the set $\Delta_P$ just by projecting the roots corresponding to $P$ to $\mathfrak{a}_P$.

Recall that the volume of the quotient space $\mathbb{Z}(\Delta_0)\backslash \mathfrak{a}_0^G$ is given by 
\[|\rm{det}(\langle\alpha_i,\alpha_j\rangle)_{1\le i,j\le 4}|^{\frac{1}{2}}=2.\]
We now normalize the Haar measure on $\mathfrak{a}_0^G=A_{G,\infty}^+\backslash A_{0,\infty}^+$ such that \[\rm{vol}(\mathbb{Z}(\Delta_0)\backslash \mathfrak{a}_0^G)=1.\]
Let $dH$ denote the corresponding  Lebesgue measure of $\mathfrak{a}_0^G$. For any test function $f\in C_c^\infty(\mathfrak{a}_0^G)$ and any element $H\in \mathfrak{a}_0^G$ expressed as $H=\sum_{i=1}^3h_ie_i$,
\begin{eqnarray}\label{fffd}
    \int_{\mathbb{R}^3}f(H)dH=\frac{1}{2}\int_{\mathbb{R}}\int_{\mathbb{R}}\int_{\mathbb{R}}f(\sum_{i=1}^3h_ie_i)dh_1\,dh_2\,dh_3,
\end{eqnarray}
it follows that the volume equals 1.

Also, if we rewrite $H=\sum_{i=1}^3a_i\alpha_i$, for Lebesgue measure $da_i,\ i=1,2,3$ on $\mathbb{R}^3$, we write (\ref{fffd}) as
\begin{eqnarray*}
    \int_{\mathbb{R}^3}f(H)dH=\frac{1}{4}\int_{\mathbb{R}}\int_{\mathbb{R}}\int_{\mathbb{R}}f(\sum_{i=1}^3a_i\alpha_i)da_1\,da_2\,da_3.
\end{eqnarray*}
An analogous normalization procedure is applied to the measures associated with any other Levi subgroup.

Finally, for parameter $T\in\mathfrak{a}_0^G$, write it as 
\[T=T_1\alpha_1+T_2\alpha_2+T_3\alpha_3.\]
\subsection{Norm}\label{2.1}
We now introduce the exact calculation of $H_P$.

Recall that we have defined the map $H_P$ from $M$ to $\mathfrak{a}_M$,
\[\langle H_P(m),\chi\rangle=\rm{log}|\chi(m)|\ (=\sum_{v\in\Sigma}\rm{log}|\chi(m_v)|_v).\] 
For $g=(g_v)\in G(\mathbb{A})$ and $P\in\mathcal{F}$, we can define the local map $H_{P,v}$, such that \[H_P(g)=\sum_{v\in\Sigma}H_{P,v}(g_v).\]
For any $n\in \mathbb{Z}_{> 0}$ and $v\in \Sigma$, $g=(g_{ij})_{n\times n}\in\rm{GL}(n,F_v)$, by Iwasawa decomposition, write \[g=nak,\quad n\in N_0,\ k\in K,\ a=\rm{diag}(a_1,a_2,...,a_n).\]
We have 
\begin{eqnarray*}
    |a_n|_v=||(g_{n1},g_{n2},...,g_{nn})||_v,\quad |a_{n-1}a_n|_v=||(\rm{det}\left(\begin{matrix}
        g_{n-1,j}&g_{n-1,k}\\
        g_{n,j}&g_{n,k}
    \end{matrix}\right),j\neq k)||_v,\ \ldots\ ,
\end{eqnarray*}
\begin{eqnarray*}
    |a_na_{n-1}\ldots a_1|_v=|\rm{det}(g_{ij})|_v.
\end{eqnarray*}
where $||(a,b,...)||_v$ is the maximum of $|a|_v,|b|_v,...$ in the non-archimedean case, and the square root of the sum of the squares of $|a|_v,|b|_v,...$ in the archimedean case.

By taking ratios, we obtain \[H_{P,v}(g)=(\rm{log}|a_1|,\rm{log}|a_2|,\ldots,\rm{log}|a_n|).\]
\section{Weighted orbital integrals and modification}\label{3}
\subsection{Weighted orbital integrals}
Now we recall the definition of weighted orbital integrals $J_M^T(\gamma,f)$ from \cite{A6} for $f\in C_c^\infty(G(F_S))$, where $S$ is a finite subset of $\Sigma$ containing $\Sigma_\infty$.

Let \( G_{\gamma,+} \) denote the centralizer of \( \gamma \) in \( G \), let \( G_\gamma \) be the Zariski-connected component of the identity in \( G_{\gamma,+} \), and let \( \mathfrak{g}_\gamma \) be its Lie algebra. 
We fix a \( G(F_S) \)-invariant measure on \( G_\gamma(F_S) \backslash G(F_S) \) and set
\[
D(\gamma) = \det\bigl(1 - \operatorname{Ad}(\sigma)\bigr)_{\mathfrak{g}/\mathfrak{g}_\gamma} \in F,
\]
where \( \sigma \) is the semisimple part of \( \gamma \). 
Assume first that $G_\gamma \subset M$, then the function \( v_M \) is left \( G_\gamma(F_S) \)-invariant, and define the weighted orbital integrals
\begin{eqnarray}\label{eqww}
    J_M^T(\gamma, f) = |D(\gamma)|_S^{1/2} \int_{G_\gamma(F_S) \backslash G(F_S)} f(g^{-1}\gamma g) \, v_M(g, T) \, dg.
\end{eqnarray}
When \( G_\gamma \not\subset M \), one uses weighted orbital integrals \( J_L^T(a\gamma, f) \) for Levi subgroups \( L \) containing \( M \) and elements \( a \in A_M(F) \) such that \( G_{a\gamma} \subset L \). 
The elements $a\gamma$ tend to $\gamma$ as $a \to 1$ in $A_M(F)$.
Arthur introduced a \( (G, M) \)-family \( \{ r_P(\lambda, \gamma, a) \mid P \in \mathcal{P}(M) \} \) so that, with the help of the corresponding function
\[
r_M^G(\gamma, a) = \lim_{\lambda \to 0} \sum_{P \in \mathcal{P}(M)} \frac{r_P(\lambda, \gamma, a)}{\theta_P(\lambda)},
\]where $\theta_P(\lambda)=\Pi_{\alpha\in\Delta_P}\langle\lambda,\alpha\rangle$.

The analogues of $r_M^G(\gamma, a)$ with \( L \) in place of \( G \), one may define
\[
J_M^T(\gamma, f) = \lim_{a \to 1} \sum_{L \in \mathcal{L}(M)} r_L^M(\gamma, a) \, J_L^T(a\gamma, f).
\]
This formula allows us to express the weighted orbital integrals of ramified orbits in terms of those of unramified orbits with a convex hull.

We shall introduce a modification to simplify the explicit formulas. 
The distribution \( J_M(\gamma, f) \) which is not dependent on $T$ is recovered by replacing all points \( T_P \) by zero. 
One then no longer has the full invariance property
\[
J_{\delta M\delta^{-1}}^T(\delta\gamma\delta^{-1}, f) = J_M^T(\gamma, f),  \qquad \delta \in G(F).
\]
\subsection{The coarse geometric expansion}
We call two elements of $G(F)$ $\mathfrak{O}$-equivalent if their semisimple parts are $G(F)$-conjugate. For any parabolic subgroup $P$, any orbit $\mathfrak{o}$ we defined in \cite{WCP}, and $f\in C_c^\infty(G(\mathbb{A})^1)$, define the partial kernel function 
\[K_{P,\mathfrak{o}}(x,y)=\sum_{\gamma\in M_P(F)\cap\mathfrak{o}}\int_{N_P(\mathbb{A})}f(x^{-1}\gamma ny)dn.\]
If we replace $K_P$ by $K_{P,\mathfrak{o}}$ in the definition (\ref{1.1}) of $J^T(F)$, we can obtain the distribution $J^T_\mathfrak{o}(f)$, and \[J^T(f)=\sum_{\mathfrak{o}\in\mathfrak{O}}J_\mathfrak{o}^T(f),\]it is called the coarse geometric expansion.

For $\rm{GL}(4)$, we have proved in \cite{WCP} that the distribution $J^T_\mathfrak{o}(f)$ is absolutely convergent. 
\subsection{The fine geometric expansion}\label{3.3}
We now state the fine geometric expansion of the distributions.

For any Levi subgroup $M$, suppose $\gamma\in M(F)$, and let $\sigma$ be the semisimple Jordan component of $\gamma$. If $\gamma'$ is another element in $M(F)$, we shall say $\gamma'$ is $(M,S)$-equivalent to $\gamma$ if there exists a $\delta\in M(F)$ with the following two properties:
\begin{itemize}
    \item $\sigma$  is also the semisimple Jordan component of $\delta^{-1}\gamma'\delta$.
    \item $\delta^{-1}\gamma$ and $\sigma^{-1}\cdot \delta^{-1}\gamma'\delta$, regarded as unipotent elements in $M_{\sigma}(F_S)$, are $M_{\sigma}(F_S)$-conjugate.
\end{itemize}
Denote by $(M(F))_{M,S}$ the set of $(M,S)$-equivalent classes in $M(F)$. Then Arthur proved in (\cite{A5}), for every compact subset $\Delta$ of $G(\mathbb{A})^1$, there exists a finite subset $S_\Delta\supset \Sigma_\infty$ of $\Sigma$ such that, for every finite subset $S$ of $\Sigma$ containing $S_\Delta$, every $M\in\mathcal{L}$ and $\gamma\in M(F)$, there exist constants $a^M(S,\gamma)$ such that 
\begin{eqnarray}
    J^T(f)=\sum_{M\in\mathcal{L}}|W_0^M||W_0^G|^{-1}\sum_{\gamma\in(M(F))_{M,S}}a^M(S,\gamma)J_M^T(\gamma,f),
\end{eqnarray}
for $f\in C_c^\infty(G(F_S)^1)$ with $\rm{supp}(f)\subset{\Delta}$ for $T=T_0$.

On the left-hand side, $f$ is extended to $G(\mathbb{A})^1$ by setting $f(g_S g^S) = f(g_S) \, \phi_{K^S}(g^S)$ for $g_S \in G(F_S)$ and $g^S \in G(\mathbb{A}^S)$, where $\phi_{K^S}$ is the characteristic function of $K^S = \prod_{v \notin S} K_v$.  
By conjugacy invariance, we may replace the sum over $\mathcal{L}$ by a sum over the set of $G(F)$-conjugacy classes of Levi subgroups $M$, in that case the factor $|W_0^M| \, |W_0^G|^{-1}$ disappears. 
\subsection{The modified weighted orbital integrals}
In this subsection, we introduce a modification to calculate the weighted orbital integrals.

Decompose the integral (\ref{eqww}) in the case $G_{a\gamma}\subset M$ by writing $g=nmk$ with $n\in N_{P_1}(F_S),m\in M_{\gamma}(F_S)\backslash M(F_S)$ and $k\in K_S=\Pi_{v\in\Sigma}K_v$. Then $\mu=m^{-1}\gamma m$ runs through the $M(F_S)$-orbit $O_{\gamma}(M(F_S))$ of $\gamma$. Then by the substitution $n^{-1}a\mu n=a\mu v$, where $v \in N_{P_1}$, we obtain the formula
\begin{eqnarray*}
    J_M^T(a\gamma,f)=|D^M(a\gamma)|_S^{\frac{1}{2}}\delta_{P_1}(a\gamma)^{\frac{1}{2}}\int_{K_S}\int_{O_\gamma(M(F_S))}\int_{P_1(F_S)}f(k^{-1}a\mu vk)v_M(n,T)dv\,d\mu\, dk,
\end{eqnarray*}where $n$ is to be regarded as a function of $v$ and $a\mu$.

Recall that Arthur defines the weighted function by the $(G,M)$-family $v_P(\lambda,n,T)$. We now give a definition of modified weighted function.

Define \[w_P(\lambda,a,\mu v,T)=v_P(\lambda,n,T)\Pi_\beta r_\beta(\lambda,\gamma,a),\]where $\beta$ runs through the reduced roots of $(N_P/(N_P\cap N_{P_1}),A_M)$ and 
\begin{eqnarray*}
    r_\beta(\lambda,\gamma,a)=|a^\beta|_S^{\lambda(\beta_\gamma)/2}|1-a^{-\beta}|_S^{\lambda(\beta_\gamma)}=|(a^\beta-1)(1-a^{-\beta})|^{\lambda(\beta_\gamma)/2}_S.
\end{eqnarray*}
Here $\beta_{\gamma} \in \mathfrak{a}_M$ is a coroot scaled by a nonnegative factor such that the limits $w_P(\lambda, 1, \mu\nu)$ exist and are nonzero for generic $\lambda$ and $\nu$. With the corresponding function
\[
w_M(1, \mu\nu, T) = \lim_{\lambda \to 0} \sum_{P \in \mathcal{P}(M)} \frac{w_P(\lambda, 1, \mu\nu, T)}{\theta_P(\lambda)},
\]
then $J_M^T(\gamma, f)$ equals
\begin{eqnarray}\label{Jm2}
    |D^M(\gamma)|_S^{\frac{1}{2}} \, \delta_{P_1}(\gamma)^{\frac{1}{2}} \int_{K_S} \int_{O_\gamma(M(F_S))} \int_{N_{P_1}(F_S)} f(k^{-1}\mu\nu k) \, w_M(1, \mu\nu, T) \, d\nu \, d\mu \, dk.
\end{eqnarray}
Arthur originally defines $r_\beta$ to be $|a^\beta - a^{-\beta}|_S^{\lambda(\beta\gamma)}$. Since this introduces additional constants into explicit formulas, we have adopted the modified definition given here. It also satisfies the Lemma 3.5 in \cite{A6}. Thus, if we define a $(G, M)$-family as the (5.1) in \cite{A6} by
\[
r_P(\lambda, \gamma, a) = \prod_\beta r_\beta(\lambda/2, \gamma, a),
\]
where $\beta$ ranges over the reduced roots of $(P, A_M)$, then the proof of the Lemma 5.3 in \cite{A6}  applies and shows that formulas (\ref{eqww}) and (\ref{Jm2}) coincide.

\section{Zeta integrals}\label{4}
This section follows the treatment of Tate \cite{T}, and we have adjusted the exposition and arguments to fit the specific aims of this paper.
\subsection{Tate integral}
For global field $F$ and its adele ring $\mathbb{A}$, $v$ is any finite place, let $\pi_v$ denote a prime element of $F_v$ and let $q_v$ denote the cardinality of the residue field of $F_v$. If $\chi_v$ is a character on $F_v^\times$, $\chi=\Pi_v\chi_v$ is a character on $F^\times\backslash \mathbb{A}^1\cong (\mathbb{R}^\times)^0F^\times\backslash \mathbb{A}^\times$, set 
\begin{eqnarray}\label{L-function}
    L_v(s,\chi_v)=L_{F,v}(s,\chi_v)=
      \begin{cases}
        (1-\chi_v(\pi_v)q_v^{-s})^{-1}, & \text{if } v<\infty\  \text{and}\ \chi_v\  \text{is\ unramified},\\
        1, & \text{if } v<\infty\ \text{and}\  \chi_v \ \text{is\  ramified},
\end{cases}
\end{eqnarray}
\begin{eqnarray*}
    L^S(s,\chi)=L_F^S(s,\chi)=\Pi_{v\notin S}L_v(s,\chi_v),\quad L(s,\chi)=L_F(s,\chi)=\Pi_{v<\infty}L_v(s,\chi_v),
\end{eqnarray*}
where $S$ denotes a finite set of places of $F$ which contains all Archimedean places.

In fact, $L^S(s,\chi)$ is absolutely convergent and holomorphic for $\rm{Re}\,s>1$ and can be meromorphically continued to the whole complex $s$-plane.

Define 
\begin{eqnarray}
    \zeta_F^S(s)=L^S(s,1_F)\quad \text{and}\quad \zeta_F(s)=L(s,1_F),
\end{eqnarray}where $1_F$ is the trivial character of $F^\times \backslash \mathbb{A}^1$. In fact  $\zeta_F^S(s)$ has a simple pole at $s=1$. And if $\chi\neq 1_F$, $L^S(s,\chi)$ is holomorphic at $s=1$.

For $\phi\in\mathcal{S}(\mathbb{A})$, define the Tate integral 
\begin{eqnarray*}
    \zeta(\phi,s,\chi)=\int_{\mathbb{A}^\times}|x|^s\chi(x)\phi(x)d^\times x=\int_{F^\times\backslash \mathbb{A}^\times}|x|^s\chi(x)\sum_{y\in F^\times}\phi(yx)d^\times x.
\end{eqnarray*}
The Tate integral is absolutely convergent and holomorphic for $\rm{Re}\, s>1$.

For fixed $S$, set $\chi_S=\Pi_{v\in S}\chi_v$. For $\phi\in\mathcal{S}(\mathbb{A})$, let \[\zeta_S(\phi_S,s,\chi_S)=\int_{F_S^\times}\phi_S(x_S)|x_S|_S^s\chi_S(x_S)d^\times x_S.\]
This integral is absolutely convergent and holomorphic for $\rm{Re}\,s>0$, and $\zeta_S(\phi_S,s,\chi_S)$ is meromorphically continued to the whole complex $s$-plane.

For $\phi\in\mathcal{S}(\mathbb{A})$, we can decompose it as \[\phi=\phi_S\times \Pi_{v\notin S}\phi_{0,v},\]where $\phi_S\in \mathcal{S}(F_S)$, and $\phi_{0,v}$ is the characteristic function of $\mathcal{O}_v$ for each $v\notin S$.

Suppose $\chi_v$ is unramified for any $v\notin S$, we have \[\zeta(\phi,s,\chi)=\zeta_S(\phi_S,s,\chi_S)\times L^S(s,\chi).\]
The residue $c_F$ of $\zeta_F(s)$ at $s=1$ satisfies $c_F=\int_{F^\times\backslash \mathbb{A}^1}d^1x$, with respect to the Tamagawa measure.

Denote 
\begin{eqnarray*}
    c_v=\begin{cases}
        (1-q_v^{-1})^{-1}, &\text{if}\  v<\infty,\\
        1, &\text{if}\  v|\infty,
    \end{cases}\quad c_S=\Pi_{v\in S}c_v,\quad \text{and}\quad c_F^S=c_S^{-1}c_F.
\end{eqnarray*}
Also, $c_Sdx_S=|x_S|_Sd^\times x_S$.

It is known that $c_F^S$ is the residue of $L^S(s,1_F)$ at $s=1$. Define the constant $\mathfrak{c}_F(S,\chi)$ by \[\mathfrak{c}_F(S,\chi)=\rm{lim}_{s\rightarrow 1}\frac{d}{ds}(s-1)L^S(s,\chi).\]
If $\chi\neq 1_F$, then $\mathfrak{c}_F(S,\chi)=L^S(1,\chi)$, define the constant $\mathfrak{c}_F'(S,\chi)$ by
\[\mathfrak{c}_F'(S,\chi)=\frac{1}{2}\rm{lim}_{s\rightarrow 1}\frac{d^2}{ds^2}(s-1)^2L^S(s,\chi),\]
similarly,\[\mathfrak{c}_F''(S,\chi)=\frac{1}{6}\rm{lim}_{s\rightarrow 1}\frac{d^3}{ds^3}(s-1)^3L^S(s,\chi).\]
We simply denote $\mathfrak{c}_F(S,1_F)$ by $\mathfrak{c}(S)$, $\mathfrak{c}_F'(S,1_F)$ by $\mathfrak{c}_F'(S)$ and $\mathfrak{c}_F''(S,1_F)$ by $\mathfrak{c}_F''(S)$.

Then we can write the Laurent expansion of $\zeta_F^S(s)$ at $s=1$ as
\[\zeta_F^S(s)=\frac{c_F^S}{s-1}+\mathfrak{c}_F(S)+\mathfrak{c}_F'(S)(s-1)+\cdots.\]

\subsection{Zeta integral}
For $n\in\mathbb{Z}_{>0}$, define $A=F^n$, and $M=(\mathbb{G}_m)^n$, so that $M(F)=(F^\times)^n$. A Haar measure $dm$ on $M(\mathbb{A})$ is defined by $dm=\Pi_{i=1}^ndm_i$, where $dm_i$ is the Haar measure on $\mathbb{A}^\times$. Let $V$ denote the vector space with the F-structure $V(F)=A$.

For $s=(s_1,s_2,...,s_n)\in \mathbb{C}^n$ and $m=(m_1,m_2,...,m_n)\in M(\mathbb{A})$, we set 
\[|m|_A^s=\Pi_{i=1}^n|m_i|_\mathbb{A}^{s_i},\quad \text{and}\quad |m|_A=\Pi_{i=1}^n|m_i|_\mathbb{A}.\]
For $m=(m_1,m_2,..,m_n)\in M$ and $x=(x_1,...,x_n)\in V$, the action of $M$ on $V$ is defined by \[m\cdot x=mx=(m_1x_1,m_2x_2,...,m_nx_n).\]
We now define a zeta integral $\zeta_A(\phi,s)$ by 
\[\zeta_A(\phi,s)=\int_{M(\mathbb{A})}|m|_A^s\phi(m)dm,\]
which is 
\[\int_{M(F)\backslash M(\mathbb{A})}|m|_A^s\sum_{x\in A^\times}\phi(mx)dm,\]for $s\in\mathbb{C}^n$, and $\phi\in\mathcal{S}(V(\mathbb{A}))$. If $\rm{Re}\  s_k>1,\,k=1,2,...,n$, then $\zeta_A(\phi,s)$ absolutely converges.

We can now define 
\[\zeta_A(\phi,s,\chi)=\int_{V_\mathbb{A}^\times}|x|^s\chi(x)\phi(x)d^\times x=\int_{A^\times\backslash V_\mathbb{A}^\times}|x|^s\chi(x)\sum_{y\in A^\times}\phi(yx)d^\times x,\]for $s\in\mathbb{C}^n$, and $\phi\in\mathcal{S}(V(\mathbb{A}))$.

We have 
\begin{eqnarray}\label{relation}
    \zeta_A(\phi,s)=\zeta_A(\phi,s,1_A).
\end{eqnarray}

\subsection{Higher rank case}
In this subsection, we introduce a basic Lemma proved by Chaudouard.

For any place $v$, let $1_v$ be the function on $F_v$ defined as follows:
\begin{itemize}
\item $1_v(x) = \exp(-\pi x^2)$ if $F_v = \mathbb{R}$;
\item $1_v(z) = \exp(-\pi |z|)$ if $F_v = \mathbb{C}$;
\item $1_v$ is the characteristic function of $\mathcal{O}_v$ if $F_v$ is non-archimedean.
\end{itemize}
Define $1_F$ to be the restricted direct product of $1_v$ with respect to the places $v$. Write $\Gamma$ as the Euler gamma function.

Let $\mathcal{S}(F)$ and $\mathcal{S}(\rm{GL}(n,\mathbb{A}))$ denote the spaces of complex-valued Schwartz--Bruhat functions on $F$ and $\rm{GL}(n,\mathbb{A})$, respectively. If $F$ is non-archimedean, these are simply the spaces of smooth compactly supported functions.

For any $f \in \mathcal{S}(F)$, set
\[
Z(f, s) = \int_{F^{\times}} f(x) |x|^s \, d^{\times}x.
\]
The integral converges and is holomorphic in $s$ for $\operatorname{Re}(s) > 0$, and admits a meromorphic continuation to $\mathbb{C}$. When $f = 1_F$, we simply write $Z(s) = Z(1_F, s)$. Explicitly, $Z(s)$ is given by 
\[
\begin{cases}
\pi^{-s/2} \Gamma(s/2), & \text{if } F = \mathbb{R},\\[4pt]
(2\pi)^{1-s} \Gamma(s), & \text{if } F = \mathbb{C},\\[4pt]
\dfrac{1}{1 - q^{-s}}, & \text{if } F \text{ is non-archimedean},
\end{cases}
\]
where $q$ is the cardinality of the residue field of $F$.

For any place $v$ of $F$, define the indentity function $1_v\in \mathcal{S}(G(F_v))$. Generally, for any finite set of places $S$, define \[1_S=\bigotimes_{v\in S}1_v\in \mathcal{S}(G(\mathbb{A}_S)),\] and \[1^S=\bigotimes_{v \notin S}1_v\in \mathcal{S}(G(\mathbb{A}^S)).\]
For any $s\in\mathbb{C}$, define \[Z_n(s)=\int_{\rm{GL}(n,\mathbb{A})}1_{F}(g)|\rm{det}\ g|^s dg.\]
One shows that this integral converges absolutely as $\rm{Re}\ s>n$. And on the open set defined by $\rm{Re}\ s>0$, we have 
\begin{eqnarray}\label{Z_n}
    Z_n(s+n)=Z_1(s+1)Z_1(s+2)\dots Z_1(s+n).
\end{eqnarray}
This yields the meromorphic continuation of $Z_n(s)$ to $\mathbb{C}$.

Also, we introduce the function\[\tilde{Z}_n(s)=(s-n)Z_n(s).\]Then $\tilde{Z}_n(s+n)$ is holomorphic for $\rm{Re}\ s>-1$.

\subsection{Truncated zeta integrals}
For a finite set of places $S$, set \[Z_{n,S}(s)=\int_{\rm{GL}(n,\mathbb{A}_S)}1_S(g)|\rm{det}\ g|^sdg,\]and \[Z_n^S(s)=\frac{Z_n(s)}{Z_{n,S}(s)}.\]Note that $Z_{n,S}(s+n)$ and $\tilde{Z}_{n}^S(s+n)=sZ_n^S(s+n)$ are holomorphic for $\rm{Re}\ s>-1$.

\begin{lemma}[\cite{C1}]
    Let $S$ be a finite set of places, and let $f\in\mathcal{S}(\rm{GL}(d,\mathbb{A}_S)^r)$. Define
    \[Z(s_1,s_2,...,s_r,f)=\int_{\rm{GL}(n,\mathbb{A}_S)^r}f(g_1,g_2,...,g_r)\Pi_{i=1}^r|\rm{det}\ g_i|^{d+s_i}dg.\]
    \begin{itemize}
        \item The integral convergent absolutely for $\rm{Re}\ s>-1$.
        \item It admits a meromorphically continuation to $\mathbb{C}^r$.
        \item Let $1=1^{\otimes r}_{\rm{GL}(d,\mathbb{A}_S)}$ and let $I\subset \{1,2,...,r\}$. Let $H_I$ be the hyperplane of $\mathbb{C}^r$ defined by the equations $s_i=0$ for all $i\in I$. Then, for all $s=(s_1,s_2,...,s_r)$ in an open neighborhood of $H_I$, one has \[Z(s_1,s_2,...,s_r,f)=Z(s_1,s_2,...,s_r,1)\frac{Z((s_i)_{i\notin I},f_I)}{Z((s_i)_{i\notin I},1_I)},\]where \[f_I((A_i)_{i\notin I})=\int_{\rm{GL}(d,\mathbb{A}_S)^I}f(g_1,g_2,...,g_r)dA_I,\]and $dA_I$ denotes the additive Haar measure on $\rm{GL}(d,\mathbb{A}_S)^I$.
    \end{itemize}
\end{lemma}

Similarly as the definition of $c_F$, $c_F^S$, $\mathfrak{c}_F(S)$, we can define $c_{2,F}$ to be the residue of $Z_2(s)$ at $s=2$, it satisfies $c_{2,F}=\int_{GL(2,F)\backslash GL(2,\mathbb{A})^1}dg$ and $c_{2,F}^S$, $\mathfrak{c}_{2,F}(S)$.

We also have $c_{2,S}d^\times g_S=|\rm{det}\,g_S|^2_Sd^+g_S$.
\section{Modified zeta integral}\label{5}
Suppose $\rm{dim}\,\mathfrak{a}_P=n$. Denote $\mathfrak{Y}_n$ by the set of subsets of $\{1,2,..,n\}$. Let $T=(T_1,T_2,...,T_n)\in \mathfrak{a}_P$, and $Y\in\mathfrak{Y}_n$. Denote the complement of $Y$ in $\{1,2,..,n\}$ by $\overline{Y}$. Set $\tau_{T,Y}(t_1,t_2,...,t_n)$ to be the characteristic function of \[\{(t_1,t_2,...,t_n)\in \mathbb{R}^n:\, t_k<-T_k,\,k\in Y\}.\]
For any vector space $V$ with $\rm{dim}\,V=n$, set \[dv_V=\Pi_{k\in V}dv_k,\quad v=(v_1,v_2,..,v_n)\in V(\mathbb{A}),\]
\[A_Y=\{(x_1,x_2,...,x_n)\in A:\,x_k=0,\,k\in \overline{Y}\},\]and 
\[V_Y(\mathbb{A})=\{(x_1,x_2,...,x_n)\in V(\mathbb{A}):\,x_k=0,\, k\in\overline{Y}\}.\]

For $v_Y=(v_1,v_2,...,v_n)\in V_Y(\mathbb{A})$, set $dv_Y=\Pi_{k\in Y}dv_k.$

We now define a modified zeta integral $\zeta_A(\phi,s,T)$ by 
\begin{eqnarray*}
    \int_{M(F)\backslash M(\mathbb{A})}dm|m|^s_A\sum_{Y\in\mathfrak{Y}_n}(-1)^{|\overline{Y}|}\int_{V_{\overline{Y}}(\mathbb{A})}\sum_{x_Y\in A_Y^\times}\phi(m(v_{\overline{Y}}+x_Y))dv_{\overline{Y}}\tau_{T,Y}(H_P(m)).
\end{eqnarray*}
Let $\psi_Y$ be a non-trivial additive character on $A_Y\backslash V_Y(\mathbb{A})$, for $\phi\in \mathcal{S}(V(\mathbb{A}))$ and $Y\in \mathfrak{Y}_n$, define 
\[\hat{\phi}^Y(x_{\overline{Y}}+x_Y)=\int_{V_Y(\mathbb{A})}\phi(x_{\overline{Y}}+y_Y)\psi_Y(x_Yy_Y)dy_Y,\]where $x_Y\in V_Y(\mathbb{A})$ and $x_{\overline{Y}}\in V_{\overline{Y}}(\mathbb{A})$.

Write $1_n$ to be the element $(1,1,...,1)\in \mathbb{C}^n$. Then by the Poisson summation formula, with respect to the self-dual Haar measure, we have
\begin{eqnarray}\label{poisson0}
    \sum_{x_Y\in A_Y}\phi(a(x_{\overline{Y}}+x_Y))=\sum_{x_Y\in A_Y}\hat{\phi}^Y(ax_{\overline{Y}}+a^{-1}x_Y)|a_Y|_A^{-1_n},
\end{eqnarray}
where $a\in M(\mathbb{A})$.
\begin{lemma}
    For any $\phi\in\mathcal{S}(V(\mathbb{A}))$, $Y\in\mathfrak{Y}_n$ and $a\in V_Y(\mathbb{A})$,
    \begin{eqnarray}\label{poisson3}
        \sum_{x_Y\in A_Y^\times}\phi(ax_Y)=\sum_{Z\subset U\subset Y}(-1)^{|Y-U|}\sum_{x_Z\in A_Z^\times}\hat{\phi}^{U}(0_{Y-Z}+a^{-1}x_Z)|a_U|_A^{-1_n}.
    \end{eqnarray}
\end{lemma}

\begin{proof}
    We prove this lemma by the induction on the dimension of $A_Y$.

    If $\rm{dim}\,A_Y=1$, we have 
    \begin{eqnarray}\label{dim1}
        \sum_{x_Y\in A_Y^\times}\phi(ax_Y)=\sum_{x_Y\in A_Y^\times}\hat{\phi}^{Y}(a^{-1}x_Y)|a_Y|_A^{-1_n}+\hat{\phi}^{Y}(0_Y)|a_Y|_A^{-1_n}-\phi(0_Y),
    \end{eqnarray}
    that is the equality (\ref{poisson3}).

    Suppose we have (\ref{poisson3}) for $\rm{dim}\,A_Y\le m-1$.

    For $\rm{dim}\, A_Y=m$, write $Y=Y_1\cup Y_2$, where $\rm{dim}\,A_{Y_1}=m-1$, $\rm{dim}\, A_{Y_2}=1$. Also, we write $\phi(x_Y)=\phi(x_{Y_1},x_{Y_2})$.
    By the induction hypothesis, 
    \[\sum_{x_{Y_1}\in A_{Y_1}^\times,x_{Y_2}\in A_{Y_2}^\times}\phi(ax_{Y_1},ax_{Y_2})\]equals
    \[\sum_{x_{Y_2}\in A_{Y_2}^\times}\sum_{Z_1\subset U_1\subset Y_1}(-1)^{|Y_1-U_1|}\sum_{x_{Z_1}\in A_{Z_1}^\times}\hat{\phi}^{U_1}(0_{Y_1-Z_1}+a^{-1}x_{Z_1},ax_{Y_2})|a_{U_1}|_A^{-1_n}.\]
    Then by (\ref{dim1}), the right hand side is \[\sum_{Z\subset U\subset Y}(-1)^{|Y-U|}\sum_{x_Z\in A_Z^\times}\hat{\phi}^{U}(0_{Y-Z}+a^{-1}x_Z)|a_U|_A^{-1_n}.\]This completes the proof.
\end{proof}
Define the sections 
\[S_{Y,T}^+=\{(t_1,t_2,...,t_n)\in ((\mathbb{R}^\times)^0)^n:\, t_k>e^{-T_k},\, k\in Y,\, t_k=1,\, k\in\overline{Y}\},\]
\[S_{Y,T}^-=\{(t_1,t_2,...,t_n)\in ((\mathbb{R}^\times)^0)^n:\, t_k<e^{-T_k},\, k\in Y,\, t_k=1,\, k\in\overline{Y}\}.\]
Let $d^\times t_Y=\Pi_{k\in Y}d^\times t_k$, since 
\[M(F)\backslash M(\mathbb{A})\cong M(F)\backslash M(\mathbb{A})^1\times A_{M,\infty}^+,\] by (\ref{poisson3}), $\zeta_A(\phi,s,T)$ equals
\begin{eqnarray*}
    &\int_{M(F)\backslash M^1(\mathbb{A})}d^1m\sum_{Y\in\mathfrak{Y}_n}(-1)^{|\overline{Y}|}\sum_{\{Z:\,Z\subset Y\}}\int_{S_{\overline{Y},T}^+\cup S_{Z,T}^-}d^\times t_{\overline{Y}\cup Z}\\
    &\int_{V_{\overline{Y}}(\mathbb{A})}\sum_{x_Z\in A_Z^\times}|t_{\overline{Y}\cup Z}|_A^s\sum_{\{U:\,Z\subset U\subset Y\}}\hat{\phi}^U(mtv_{\overline{Y}}+m^{-1}t^{-1}x_Z+0_{Y-Z})dv_{\overline{Y}}\\
    &\cdot (-1)^{|Y-U|}|t_Z|_A^{-1_n}\Pi_{k\in Y-Z}\frac{e^{-s_kT_k}}{s_k}\Pi_{j\in U-Z}\frac{e^{-(s_j-1)T_j}}{s_j-1}.
\end{eqnarray*}
For fixed $Y\in\mathfrak{Y}_n$, write $\phi_Y(x_Y)=\int_{V_{\overline{Y}}(\mathbb{A})}\phi(x_Y+x_{\overline{Y}})dv_{\overline{Y}}$, $s_Y=(s_k)_{k\in Y}$, and $T_Y=(T_k)_{k\in Y}$.

Define $\zeta_{A_Y}(\phi_Y,s_Y,T_Y)$ by
\begin{eqnarray}\label{zetaa}
    &\int_{M(F)\backslash M^1(\mathbb{A})}d^1m\sum_{\{Z:\,Z\subset Y\}}\int_{S_{Z,T}^-}d^\times t_{Z}\int_{V_{\overline{Y}}(\mathbb{A})}\sum_{x_Z\in A_Z^\times}|t_{Z}|_A^s\sum_{\{U:\,Z\subset U\subset Y\}}\\\notag
    &\hat{\phi}^U(m^{-1}t^{-1}x_Z+0_{Y-Z})dv_{\overline{Y}}\cdot (-1)^{|Y-U|}|t_Z|_A^{-1_n}\Pi_{k\in Y-Z}\frac{e^{-s_kT_k}}{s_k}\Pi_{j\in U-Z}\frac{e^{-(s_j-1)T_j}}{s_j-1}.
\end{eqnarray}
After replacing $mtv_{\overline{T}}$ by $v_{\overline{T}}'$ in $\zeta_A(\phi,s,T)$, the limit  
\[\rm{lim}_{s\rightarrow 1_n}\Pi_{i=1}^n\frac{\partial}{\partial s_i}\Pi_{j=1}^n (s_j-1) \zeta_A(\phi,s,T)\]equals
\[\rm{lim}_{s\rightarrow 1_n}\Pi_{i=1}^n\frac{\partial}{\partial s_i}\Pi_{j=1}^n (s_j-1) \sum_{Y\in\mathfrak{Y}_n}\Pi_{l\in \overline{Y}}|T_l|\zeta_{A_{Y}}(\phi_Y,s_{|Y|},T_Y).\]

Recall that $\zeta_A(\phi,s)$ equals
\begin{eqnarray*}
    \int_{M(F)\backslash M(\mathbb{A})}|m|_A^s\sum_{x\in A^\times}\phi(mx)dm.
\end{eqnarray*} 
By (\ref{poisson0}), we can write $\zeta_A(\phi,s)$ as 
\begin{eqnarray*}
    \int_{M(F)\backslash M(\mathbb{A})}|m|_A^s\sum_{Z\subset U\in\mathfrak{Y}_n}(-1)^{|\overline{U}|}\sum_{x_Z\in A_Z^\times}\hat{\phi}^U(0_{\overline{Z}}+m^{-1}x_Z)|m|^{-1_n}_A d^1hdm.
\end{eqnarray*}
We have \[|\hat{\phi}^Y(m(v_{\overline{Y}}+x_Y))\tau_{Y,T}(H_P(m))|\le|\hat{\phi}^Y(m(v_{\overline{Y}}+x_Y))|,\]
since $\tau_{Y,T}(H_P(m))\le 1$.

Define $\zeta_A^T(\phi,s)$ by 
\[\int_{M(F)\backslash M(\mathbb{A})}|m|_A^s\sum_{Z\subset U\in\mathfrak{Y}_n}(-1)^{|\overline{U}|}\sum_{x_Z\in A_Z^\times}\hat{\phi}^U(0_{\overline{Z}}+m^{-1}x_Z)\tau_{Y,T}(H_P(m))|m|^{-1_n}_A d^1hdm.\]
Thus $\zeta_A^T(\phi,s)$ equals 
\begin{eqnarray}\label{zeta222}
    &\int_{M(F)\backslash M(\mathbb{A})}dm\sum_{Z\subset U\in\mathfrak{Y}_n}\int_{S_{Z,T}^-}d^\times t_{Z}\\\notag
    &\sum_{x_Z\in A_Z^\times}|t_{Z}|_A^s\hat{\phi}^U(m^{-1}t^{-1}x_Z+0_{\overline{Z}})\cdot (-1)^{|\overline{U}|}|t_Z|_A^{-1_n}\Pi_{k\in \overline{U}}\frac{e^{-s_kT_k}}{s_k}\Pi_{j\in U-Z}\frac{e^{-(s_j-1)T_j}}{s_j-1}.
\end{eqnarray}
Compare the formula of $\zeta_A(\phi,s,T)$ and $\zeta^T_A(\phi,s)$, we have 
\begin{theorem}
    For $T=T_0=(0,0,...,0)\in \mathfrak{a}_P$, $\phi\in\mathcal{S}(V(\mathbb{A}))$, and $s\in\mathbb{C}^n$,
    \[\rm{lim}_{s\rightarrow 1_n}\Pi_{i=1}^n\frac{\partial}{\partial s_i}\Pi_{j=1}^n (s_j-1) \zeta_A(\phi,s,T_0)=\rm{lim}_{s\rightarrow 1_n}\Pi_{k=1}^n\frac{\partial}{\partial s_k}\Pi_{j=1}^n(s_j-1)\zeta_A^{T_0}(\phi,s). \]
\end{theorem}
If we write $\zeta_A(\phi,s,1_A)$ as the product \[\int_{M(F_S)}|m_S|_A^s\phi_S(m)dm_S\times L^S_{A}(\phi,s),\]where $L^S_{A}(\phi,s)$ is defined similarly as (\ref{L-function}), according to (\ref{relation}), we have 
\begin{corollary}
    For $T=T_0=(0,0,...,0)\in \mathfrak{a}_P$, $\phi\in\mathcal{S}(V(\mathbb{A}))$, and $s\in\mathbb{C}^n$, the limit 
    \[\rm{lim}_{s\rightarrow 1_n}\Pi_{i=1}^n\frac{\partial}{\partial s_i}\Pi_{j=1}^n (s_j-1) \zeta_A(\phi,s,T_0)\] equals
    \begin{eqnarray*}
    \sum_{Y\in\mathfrak{Y}_n}(c_F^S)^{|Y|}\mathfrak{c}_F(S)^{|\overline{Y}|}\int_{M(F_S)}\phi_S(m_S)|m_S|_A\cdot(\Pi_{k\in Y}\rm{log}\,|m_k|_{A})dm_S.
    \end{eqnarray*}
\end{corollary}
\begin{proof}
    According to (\ref{relation}), we consider $\zeta_A(\phi,s,1_A)$.

    Write \[f(s)=\int_{M(F)\backslash M(F_S)}|m_S|_A^s\phi_S(m)dm_S,\]then \[\zeta_A(\phi,s,1_A)=f(s)\times L^S_{A}(\phi,s).\]Since $f(s)$ is holomorphic at $s=1$, write
    $\Pi_{j=1}^n (s_j-1) \zeta_A(\phi,s,1_A)$ as
    \[f(s)\times \Pi_{j=1}^n(s_j-1)L^S_{A}(\phi,s).\]
    For any $1\le i\le n$, $\frac{\partial}{\partial s_i}\zeta_A(\phi,s,1_A)$ equals
    \begin{eqnarray*}
        \frac{\partial}{\partial s_i}f(s)\times \Pi_{j=1}^n(s_j-1)L^S_{A}(\phi,s)+f(s)\times \frac{\partial}{\partial s_i}\Pi_{j=1}^n(s_j-1)L^S_{A}(\phi,s).
    \end{eqnarray*}
    Since \[c_F^S=\rm{lim}_{s_i\rightarrow 1}\Pi_{j=1}^n(s_j-1)L^S_{A}(\phi,s),\quad \mathfrak{c}_F(S)=\rm{lim}_{s_i\rightarrow 1}\frac{\partial}{\partial s_i}\Pi_{j=1}^n(s_j-1)L^S_{A}(\phi,s),\] by the calculations through differentiation, we can obtain the desired result.
\end{proof}

\section{Mean value formula for $\rm{GL}(3)$}\label{6}
In this section, we consider $G=\rm{GL}(3)$. Let $V$ denotes the vector space with the $F$-structure $V(F)=F\oplus F\oplus F$. Let $dg$ denotes a Haar measure on $G(\mathbb{A})$. Denote by $dx_i,\, i=1,2,3$ the Haar measure on $\mathbb{A}$. Thus $dx=dx_1dx_2dx_3$ is a Haar measure on $V(\mathbb{A})$.

Define a right action $G$ on $V$ by 
\begin{eqnarray*}
   x\cdot g=(ax_1+dx_2+gx_3,bx_1+ex_2+hx_3,cx_1+fx_2+ix_3)
\end{eqnarray*}for 
\[ g=\left(\begin{matrix}
    a&b&c\\
    d&e&f\\
    g&h&i
   \end{matrix}\right)\in G,\quad x=(x_1,x_2,x_3)\in V.\]
Define $V^0$ to be the scheme-theoretic complement of $\{0\}$ in $V$.

Recall the mean value formula introduced in \cite{Si},
\begin{eqnarray}\label{mean1}
    \int_{G(F)\backslash G(\mathbb{A})}\sum_{x\in V^0(F)}f(x\cdot g)dg=\rm{vol}_G\int_{V(\mathbb{A})}f(x)dx,
\end{eqnarray}
where $f\in C_c^\infty(V(\mathbb{A}))$. And the complement of $V^0(\mathbb{A})$ in $V(\mathbb{A})$ has measure zero.

Set 
\begin{eqnarray*}
    M_0=\left\{\left(\begin{matrix}
        *&&\\
        &*&\\
        &&*
    \end{matrix}\right)\in G\right\},\quad N_0=\left\{\left(\begin{matrix}
        1&*&*\\
        &1&*\\
        &&1
    \end{matrix}\right)\in G\right\}.
\end{eqnarray*}
Fix the Haar measure on $A_{M,\infty}^+$ as $dt_1dt_2dt_3$ for $(t_1,t_2,t_3)\in A_{M,\infty}^+$, where $dt_i,\,i=1,2,3$ is the Lebesgue measure on $\mathbb{R}$.

Set \[K_3=\Pi_v K_{3,v},\quad K_{3,v}=\begin{cases}
    U(3), & \text{if}\  v\in \Sigma_\mathbb{C},\\
    O(3), & \text{if}\  v\in\Sigma_{\mathbb{R}},\\
    GL(3,\mathcal{O}_v), &\text{if}\  v\,<\,\infty.
\end{cases}\]
We have $K_3$ is equal to the stablizier of the height function $||\cdot||$ on $V^0(\mathbb{A})$ in $G(\mathbb{A})$. Fix the Haar measure $dk\,(\text{resp}.\  dn)$ on $K_3\,(\text{resp}.\  N(\mathbb{A}))$ normalized by \[\int_Kdk=1,\quad (\text{resp}.\quad \int_{N(F)\backslash N(\mathbb{A})}dn=1).\]
Since $V^0(F)$ is a single $G(F)$-orbit and $N(F)$ is the stabilizer of $(0,0,1)$, the formula (\ref{mean1}) is equivalent to 
\begin{eqnarray}\label{mean2}
    \int_{K_3}\int_{\mathbb{A}^\times}f((0,0,a)\cdot k)|a|^3da\,dk=\rm{vol}_G\int_{V(\mathbb{A})}f(x)dx.
\end{eqnarray}
If $d^\times z$ is the Haar measure on $\mathbb{A}^\times$, we have $da=c_F^{-1}\rm{vol}_Md^\times z$.
\section{ramified terms}\label{7}
In this section, we mainly discuss the case $\mathfrak{o}_{1111}^4$, and in this case, we rewrite the distribution $J_\mathfrak{o}(f)$.

Suppose first that $\mathfrak{o}$ consists entirely of unipotent elements. Then $\mathfrak{o} \subset U_G(F)$, where $U_G$ is the closed variety of unipotent elements in $G$. This class in $\mathfrak{o}_{1111}^4$ is the furthest from being unramified, and consequently the most difficult to handle.

Recall that for any orbit $\mathfrak{o}$, $f\in C_c^\infty(G(\mathbb{A})^1)$, we have defined $I^\mathfrak{o}(f,x,T)$(\cite{WCP}), which equals 
\begin{eqnarray*}
    \sum_P(-1)^{\rm{dim}\,Z\backslash A}\sum_{\delta\in P(F)\backslash G(F)}\sum_{\gamma\in M^\mathfrak{o}}\int_{N(\mathbb{A})}f(x^{-1}\delta^{-1}\gamma n\delta x)dn\cdot \hat{\tau}_P(H_0(\delta x)-T).
\end{eqnarray*}
We have proved in \cite{WCP} this function of $x$ is absolutely integrable over $Z_\infty^+G(F)\backslash G(\mathbb{A})$.

In the case of $\mathfrak{o} \subset U_G(F)$, we shall write
\begin{eqnarray}\label{ogamma}
    J_\mathfrak{o}^T(f)=\sum_{M\in\mathscr{L}(M_0)}|W_{0}^M||W_0^G|^{-1}\sum_{u\in (\mathcal{U}_M(F))_{M,S}}a^M(S,u)J_M^T(u,f).
\end{eqnarray}
Set 
\begin{eqnarray*}
    u(x_1,x_2,x_3,x_4,x_5,x_6)=\left(\begin{matrix}
        1&x_1&x_2&x_3\\
        &1&x_4&x_5\\
        &&1&x_6\\
        &&&1
    \end{matrix}\right).
\end{eqnarray*}
For each $z\in Z(F)$, denote by $\mathfrak{o}_{1111,z}^4$ the elements in $\mathfrak{o}_{1111}^4$ whose eignvalue is $z$.

We divide $\mathfrak{o}_{1111,z}^4$ into five classes by the Jordan normal form and write \[\mathfrak{o}_{1111,z}^4=\mathfrak{o}_{4,z}^4\cup\mathfrak{o}^4_{31,z}\cup\mathfrak{o}^4_{22,z}\cup\mathfrak{o}^4_{211,z}\cup\mathfrak{o}^4_{1,z},\]where 
\begin{eqnarray*}
    \mathfrak{o}_{4,z}^4=\{z\}_G,\quad  \mathfrak{o}_{31,z}^4=\{zu(0,0,1,0,0,0)\}_G,\quad  \mathfrak{o}_{22,z}^4=\{zu(0,1,0,0,1,0)\}_G,
\end{eqnarray*}
\begin{eqnarray*}
    \mathfrak{o}_{211,z}^4=\{zu(0,1,0,0,0,1)\}_G,\quad  \mathfrak{o}_{1,z}^4=\{zu(1,0,0,1,0,1)\}_G.
\end{eqnarray*}
For briefy, we omit the superscript $4$ if it does not cause confusion.

For any parabolic subgroup $P$, any set $O$ consists of elements in $M$, define \[K_{P,O}(g,h)=\sum_{\gamma\in M\cap O}\int_{N(\mathbb{A})}f(g^{-1}\gamma nh)dn,\]where $g,h\in G(\mathbb{A})$.

For any element $g\in G$, by the Langlands decomposition, for any parabolic subgroup $P=MN$, we can uniquely write \[g=m_gn_gk_g,\quad m_g\in M,\ n_g\in N,\ k\in K.\] Define a map $\rm{La}_M(g)$ from $G$ to any Levi subgroup $M$, such that $\rm{La}_M(g)=m_g$.

Write \[J_{\mathfrak{o}_{4,z}}(f)=\int_{G(F)\backslash G(\mathbb{A})^1}K_{G,\mathfrak{o}_{4,z}}(g,g)dg.\]
It is clear that \[J_{\mathfrak{o}_{4,z}}(f)=\rm{vol}_G f(z).\]
Define $f_{K_S}(g)$ to be \[\int_{K_S}f(k_S^{-1}gk_S)dk_S,\quad g\in G(\mathbb{A}).\]
Define 
\begin{eqnarray*}
    &J_{\mathfrak{o}_{31,z}}^T(f)=\int_{G(F)\backslash G(\mathbb{A})^1}[K_{G,\mathfrak{o}_{31,z}}(x,x)-\sum_{\delta\in P_{31}(F)\backslash G(F)}K_{P_{31,\{z\}}}(\delta x,\delta x)\hat{\tau}_{P_{31}}(H_0(\delta x)-T)\\
    &-\sum_{\delta\in P_{13}(F)\backslash G(F)}K_{P_{13,\{z\}}}(\delta x,\delta x)\hat{\tau}_{P_{13}}(H_0(\delta x)-T)]dx,
\end{eqnarray*}
\begin{eqnarray*}
    J_{\mathfrak{o}_{22,z}}^T(f)=\int_{G(F)\backslash G(\mathbb{A})^1}K_{G,\mathfrak{o}_{22,z}}(x,x)dx,
\end{eqnarray*}
\begin{eqnarray*}
    &J_{\mathfrak{o}_{211,z}}^T(f)=\int_{G(F)\backslash G(\mathbb{A})^1}[K_{G,\mathfrak{o}_{211,z}}(x,x)-\sum_{\delta\in P_{22}(F)\backslash G(F)}K_{P_{22},\{zu(0,0,0,0,0,1)\}}(\delta x,\delta x)\\ &\cdot\hat{\tau}_{P_{22}}(H_0(\delta x)-T)-\sum_{\delta\in P_{31}(F)\backslash G(F)}K_{P_{31},\{zu(0,1,0,0,0,0)\}}(\delta x,\delta x)\hat{\tau}_{P_{31}}(H_0(\delta x)-T)\\
    &+\sum_{\delta \in P_{211}(F)\backslash G(F)}K_{P_{211},\{z\}}(\delta x,\delta x)\hat{\tau}_{P_{211}}(H_0(\delta x)-T)]dx,
\end{eqnarray*}and
\begin{eqnarray*}
    &J_{\mathfrak{o}_{1,z}}^T(f)=\int_{G(F)\backslash G(\mathbb{A})^1}\{K_{G,\mathfrak{o}_{1,z}}(x,x)\\
    &+\sum_{P\neq G}(-1)^{\rm{dim}\,Z\backslash A}\sum_{\delta\in P(F)\backslash G(F)}K_{P,\{\rm{La}_M(zu(1,0,0,1,0,1))\}_M}(\delta x,\delta x)\hat{\tau}_P(H_0(\delta x)-T)\}dx.
\end{eqnarray*}
According to \cite{WCP}, the other terms in \[\int_{Z_\infty^+G(F)\backslash G(\mathbb{A})}\sum_{P}(-1)^{\rm{dim}\, Z\backslash A}\sum_{\delta \in P(F)\backslash G(F)}K_{P}(\delta x,\delta x)\hat{\tau}_{P}(H_0(\delta x)-T)dx\]are equal to zero if we take $T=0$.
\subsection{$\mathfrak{o}_{31,z}$ in $J_{\mathfrak{o}_{1111,z}}^T(f)$}
In this subsection, we shall calculate the elements in $\mathfrak{o}_{31,z}$ in $J_{\mathfrak{o}_{1111,z}}^T(f)$.

\begin{theorem}\label{left31}
    For $f\in C_c^\infty(G(\mathbb{A})^1)$, $z\in Z(F)$, $J_{\mathfrak{o}_{31,z}}(f)$ equals the sum of 
    \begin{eqnarray*}
        \rm{vol}_{M_{31}}\int_{F_S^3}f_{K_S}(zu(0,0,n_{14},0,n_{24},n_{34}))(\rm{log}||(n_{14},n_{24},n_{34})||)dn_{14}dn_{24}dn_{34}
    \end{eqnarray*}and \[\rm{vol}_{M_{31}}\cdot \frac{\frac{d}{ds}\zeta_F^S(s)|_{s=3}}{\zeta_F^S(3)}\cdot\int_{N_{31}(\mathbb{A})}f_{K_S}(zn)dn.\]
\end{theorem}
\begin{proof}
    Write $u_z=zu(0,0,1,0,0,0)$, then \[K_{G,\mathfrak{o}_{31,z}}(x,x)=\sum_{\gamma\in{G^{\mathfrak{o}_{31,z}}}}f(x^{-1}\gamma x)=\sum_{\delta\in G_{u_z}(F)\backslash G(F)}f(x^{-1}\delta^{-1}u_z\delta x),\]where 
    \begin{eqnarray*}
        G_{u_z}(F)=\left\{\left(\begin{matrix}
            a&&\\
            &M_{22}&\\
            &&a
        \end{matrix}\right)N_{121}(F),\quad a\in F^\times,\ M_{22}\in \rm{GL}(2,F)\right\}.
    \end{eqnarray*}
    Decompose $J_{\mathfrak{o}_{31,z}}^T(f)$ as the sum 
    \begin{eqnarray}\label{311}
        \notag \int_{G(F)\backslash G(\mathbb{A})^1}\sum_{\delta\in G_{u_z}(F)\backslash G(F)}f(x^{-1}\delta^{-1}u_z\delta x)(1-\hat{\tau}_{P_{31}}(H_0(\delta x)-T)-\hat{\tau}_{P_{13}}(H_0(\delta x)-T))dx\\
    \end{eqnarray}
    \begin{eqnarray}\label{312}
        &\notag+\int_{G(F)\backslash G(\mathbb{A})^1}\{\sum_{\delta\in G_{u_z}(F)\backslash G(F)}f(x^{-1}\delta^{-1}u_z\delta x)-\sum_{P_{31}(F)\backslash G(F)}\\
        &\int_{N_{31}(\mathbb{A})}f(x^{-1}\delta^{-1}zn\delta x)dn\}\cdot\hat{\tau}_{P_{31}}(H_0(\delta x)-T)dx
    \end{eqnarray}
    \begin{eqnarray}\label{313}
        &\notag+\int_{G(F)\backslash G(\mathbb{A})^1}\{\sum_{\delta\in G_{u_z}(F)\backslash G(F)}f(x^{-1}\delta^{-1}u_z\delta x)-\sum_{P_{13}(F)\backslash G(F)}\\
        &\int_{N_{13}(\mathbb{A})}f(x^{-1}\delta^{-1}zn\delta x)dn\}\cdot\hat{\tau}_{P_{13}}(H_0(\delta x)-T)dx.
    \end{eqnarray}
    By the Langlands decomposition, \[x=nmak,\quad n\in N_{31}(\mathbb{A}),\ m\in M_{31}(\mathbb{A})^1,\ a\in A_{31,\infty}^{G,+},\ k\in K,\]since $G_{u_z}\subset P_{121}\subset P_{31}$, 
    \[\int_{G(F)\backslash G(\mathbb{A})^1}\sum_{\delta\in G_{u_z}(F)\backslash G(F)}f(x^{-1}\delta^{-1}u_z\delta x)\]
    becomes
    \begin{eqnarray*}
        \int_{A_{31,\infty}^{G,+}}\int_{P_{31}(F)\backslash P(\mathbb{A})^1}\sum_{\delta\in G_{u_z}(F)\backslash P_{31}(F)}f_K(a^{-1}p^{-1}\delta^{-1}u_z\delta pa)dpda.
    \end{eqnarray*}
    For fixed $a\in A_{31,\infty}^{G,+},\ p\in P_{31}(F)\backslash P_{31}(\mathbb{A})^1$, \[\sum_{\delta\in G_{u_z}(F)\backslash P_{31}(F)}f_K(a^{-1}p^{-1}\delta^{-1}u_z\delta pa)=\sum_{v\in N_{31}(F),\ v\neq e}f_K(za^{-1}p^{-1}vpa).\]
    By (\ref{mean1}), 
    \begin{eqnarray*}
        \int_{A_{31,\infty}^{G,+}}\int_{P_{31}(F)\backslash P(\mathbb{A})^1}\sum_{\delta\in G_{u_z}(F)\backslash P_{31}(F)}f_K(a^{-1}p^{-1}\delta^{-1}u_z\delta pa)dpda,
    \end{eqnarray*}
    equals
    \begin{eqnarray*}
        \rm{vol}_{M_{31}}\int_{N_{31}(\mathbb{A})}f_K(za^{-1}na)dn,
    \end{eqnarray*}
    which is the same as \[\int_{P_{31}(F)\backslash P_{31}(\mathbb{A})^1}\int_{N_{31}(\mathbb{A})}f_K(a^{-1}p^{-1}znpa)dndp.\]
    Thus (\ref{312}) equals zero. Similarly, (\ref{313}) equals zero.

    As for the term (\ref{311}), replace $G(F)$ in (\ref{311}) by $P_{121}(F)$, then \[\sum_{\delta\in G_{u_z}(F)\backslash G(F)}f(x^{-1}\delta^{-1}u_z\delta x)=\sum_{v\in F^\times}f(x^{-1}zu(0,0,v,0,0,0)x).\]
    Suppose the Iwasawa decomposition\[x=n\rm{diag}(a,1,1,d)k,\quad n\in N_{121},\ a,d\in F^\times\backslash \mathbb{A}^\times,\ k\in K,\]the modular function $e^{-\langle 2\rho_{P_{121}},H_0(x)\rangle}$ equals $|d/a|^3$.

    Denote $\chi$ the characteristic function in $\mathfrak{a}_{P_{121}}^{G}$, (\ref{311}) equals 
    \begin{eqnarray*}
        \frac{1}{2}\int_{F^\times\backslash \mathbb{A}^\times}\int_{F^\times \backslash \mathbb{A}^\times}\sum_{v\in F^\times}f_K(zu(0,0,v\cdot d/a,0,0,0))|d/a|^3\\
        \cdot(1-\chi(\rm{log}|a/d^3|>4T_3)-\chi(\rm{log}|a^3/d|>4T_1))d^\times ad^\times d.
    \end{eqnarray*}
    Take a change of variables, it becomes
    \begin{eqnarray*}
        \frac{1}{2}\int_{F^\times\backslash \mathbb{A}^\times}\int_{F^\times \backslash \mathbb{A}^\times}\sum_{v\in F^\times}f_K(zu(0,0,dv,0,0,0))|d|^3\\
        \cdot(1-\chi(\rm{log}|1/d^3a^2|>4T_3)-\chi(\rm{log}|a^2/d|>4T_1))d^\times ad^\times d,
    \end{eqnarray*}
    which is 
    \begin{eqnarray*}
        \int_{F^\times \backslash \mathbb{A}^\times}\sum_{v\in F^\times}f_K(zu(0,0,dv,0,0,0))|d|^3(\rm{log}|d|+T_1+T_3)d^\times d.
    \end{eqnarray*}
    By (\ref{mean2}), (\ref{311}) equals 
    \[\rm{vol}_{M_{31}}\int_{\mathbb{A}^3}f_K(zu(0,0,n_{14},0,n_{24},n_{34}))(\rm{log}||(n_{14},n_{24},n_{34})||)dn_{14}dn_{24}dn_{34},\] if we take $T=T_0$.

    Write \[\rm{log}||(n_{14},n_{24},n_{34})||=\rm{log}||(n_{14},n_{24},n_{34})||_S+\rm{log}||(n_{14},n_{24},n_{34})||_{\overline{S}},\]where $\overline{S}$ is the complement of $S$.

    For each $v\notin S$, \[\int_{N_{31}(\mathcal{O}_v)}\rm{log}||(n_{14,v},n_{24,v},n_{34,v})||_vdn_{14,v}dn_{24,v}dn_{34,v}=\frac{q_v^{-3}\rm{log}(q_v^{-1})}{1-q_v^{-3}},\]and\[\rm{vol}_{M_{31}}\int_{N_{31}(\mathbb{A})}f_{K_S}(zn)dn.\]
    Also, we can easily obatin that \[\frac{\frac{d}{ds}\zeta_F^S(s)|_{s=3}}{\zeta_F^S(3)}=\frac{q_v^{-3}\rm{log}(q_v^{-1})}{1-q_v^{-3}},\]if we consider $g(s)=\rm{log}(\zeta_F^S(s))$.
\end{proof}
\subsection{$\mathfrak{o}_{22,z}$ in $J_{\mathfrak{o}_{1111,z}}^T(f)$}In this subsection, we shall calculate the elements in $\mathfrak{o}_{22,z}$ in $J_{\mathfrak{o}_{1111,z}}^T(f)$.

Denote $u_z=zu(0,1,0,0,1,0)$, for $z\in Z(F)$.

Recall that \[J_{\mathfrak{o}_{22,z}}(f)=\int_{G(F)\backslash G(\mathbb{A})^1}\sum_{\delta\in G_{u_z}(F)\backslash G(F)}f(x^{-1}\delta^{-1}u_z\delta x)dx.\]
Since \begin{eqnarray*}
    G_{u_z}=\left(\begin{matrix}
        A&N\\
        &A
    \end{matrix}\right)\subset P_{22},\quad A,\ N\in \rm{GL}(2),
\end{eqnarray*}
we can write $J_{\mathfrak{o}_{22,z}}(f)$ as \[\int_{P_{22}(F)\backslash G(\mathbb{A})^1}\sum_{\delta\in G_{u_z}(F)\backslash P_{22}(F)}f(x^{-1}\delta^{-1}u_z\delta x)dx.\]
For fixed $x\in P_{22}(F)\backslash G(\mathbb{A})^1$, since $N_{22}$ is abelian, \[\sum_{\delta\in G_{u_z}(F)\backslash P_{22}(F)}f(x^{-1}\delta u_z\delta x)=\sum_{A,B\in \rm{GL}(2,F)}f(x^{-1}u^{\rm{GL}(2)}(AB)x),\]where \begin{eqnarray*}
    u^{\rm{GL}(2)}(M)=\left(\begin{matrix}
        I_2&M\\
        &I_2
    \end{matrix}\right)\in P_{22},\quad M\in \rm{GL}(2).
\end{eqnarray*}
For any 
\begin{eqnarray*}
    m=\left(\begin{matrix}
        M_1&\\
        &M_2^{-1}
    \end{matrix}\right)\in M_{22}(\mathbb{A}),\quad M_1,M_2\in \rm{GL}(2,\mathbb{A}),
\end{eqnarray*}
the modular function $e^{-\langle 2\rho_{P_{22}},H_0(m)\rangle}$ equals $|\rm{det}(M_1M_2)|^2$.

By the Langlands decomposition,\[x=mn\rm{diag}(a,a,1,1)k\in P_{22}(F)\backslash G(\mathbb{A})^1.\]
$J_{\mathfrak{o}_{22,z}}(f)$ equals 
\begin{eqnarray*}
    \int_{F^\times\backslash \mathbb{A}^\times}\int_{(Z(\mathbb{A})\rm{GL}(2,F)\backslash \rm{GL}(2,\mathbb{A}))^2}f_K(zu^{\rm{GL}(2)}(aM_1M_2))|\rm{det}(aM_1M_2)|^2dadM_1dM_2,
\end{eqnarray*}
If we take a change of variables, it equals \[\rm{vol}_{\rm{GL}(2)}\int_{\rm{GL}(2,\mathbb{A})}f_K(zu^{\rm{GL}(2)}(M_1))|\rm{det}(M_1)|^2dM_1.\]
For $s\in\mathbb{C}$, define $J_{\mathfrak{o}_{22,z}}(f,s)$ by \[\rm{vol}_{\rm{GL}(2)}\int_{\rm{GL}(2,\mathbb{A})}f_K(zu^{\rm{GL}(2)}(M_1))|\rm{det}(M_1)|^{s}dM_1.\]By (\ref{Z_n}), $s=0$ is the pole.

For fixed finite set $S$ of the places, $J_{\mathfrak{o}_{22}}(f,s+2)$ can be written as the product
\begin{eqnarray*}
    \rm{vol}_{\rm{GL}(2)}\int_{\rm{GL}(2,\mathbb{A}_S)}f_{K_S}(zu^{\rm{GL}(2)}(M_S))|\rm{det}(M_S)|^{s+2}dM_S\times Z^S_2(s+2),
\end{eqnarray*}
it is clear that the former one is finite if we take $s=0$.
Consider
\begin{eqnarray*}
    \rm{lim}_{s\rightarrow 0}\frac{d}{ds}(s\cdot J_{\mathfrak{o}_{22,z}}(f,s+2)),
\end{eqnarray*}which is 
\begin{eqnarray*}
    \rm{lim}_{s\rightarrow 0}\frac{d}{ds}\{ \rm{vol}_{\rm{GL}(2)}\int_{\rm{GL}(2,\mathbb{A}_S)}f_{K_S}(zu^{\rm{GL}(2)}(M_S))|\rm{det}(M_S)|^{s+2}dM_S\times (s\cdot Z^S_2(s+2))\}.
\end{eqnarray*}
We have
\begin{theorem}
    For $f\in C_c^\infty(G(\mathbb{A})^1)$, $z\in Z(F)$, $J_{\mathfrak{o}_{22,z}}(f)$ equals the sum of 
    \begin{eqnarray*}
    \rm{vol}_{\rm{GL}(2)}\int_{\rm{GL}(2,\mathbb{A}_S)}f_{K_S}(zu^{\rm{GL}(2)}(M_S))|\rm{det}(M_S)|^2\rm{log}|\rm{det}M_S|dM_S\cdot c^S_{2,F}
\end{eqnarray*}and \[\rm{vol}_{\rm{GL}(2)}\mathfrak{c}_{2,F}(S)\int_{\rm{GL}(2,\mathbb{A}_S)}f_{K_S}(zu^{\rm{GL}(2)}(M_S))|\rm{det}(M_S)|^{2}dM_S.\]
\end{theorem}
\subsection{$\mathfrak{o}_{211,z}$ in $J_{\mathfrak{o}_{1111,z}}^T(f)$}In this subsection, we shall calculate the elements in $\mathfrak{o}_{211,z}$ in $J_{\mathfrak{o}_{1111,z}}^T(f)$.

Let $V$ denote the vector space with the $\mathbb{F}$-structure $V[\mathbb{F}]=F\oplus F$. Set $V^0$ the scheme-theoretic complement of $(0,0)$ in $V$. Denote $u_z=zu(0,1,0,0,0,1)$, for $z\in Z(F)$.

Recall that $J_{\mathfrak{o}_{211,z}}^T(f)$ equals
\begin{eqnarray*}
    &\int_{G(F)\backslash G(\mathbb{A})^1}[K_{G,\mathfrak{o}_{211,z}}(x,x)-\sum_{\delta\in P_{22}(F)\backslash G(F)}K_{P_{22},\{zu(0,0,0,0,0,1)\}}(\delta x,\delta x)\\
    &\cdot \hat{\tau}_{P_{22}}(H_0(\delta x)-T)-\sum_{\delta\in P_{31}(F)\backslash G(F)}K_{P_{31},\{zu(0,1,0,0,0,0)\}}(\delta x,\delta x)\hat{\tau}_{P_{31}}(H_0(\delta x)-T)\\
    &+\sum_{\delta \in P_{211}(F)\backslash G(F)}K_{P_{211},\{z\}}(\delta x,\delta x)\hat{\tau}_{P_{211}}(H_0(\delta x)-T)]dx.
\end{eqnarray*}
Firstly, we need to prove
\begin{lemma}
    The integral 
    \begin{eqnarray}\label{+211}
        \int_{G(F)\backslash G(\mathbb{A})^1}\sum_{\delta\in P_{31}(F)\backslash G(F)}K_{P_{31},\{zu(0,1,0,0,0,0)\}}(\delta x,\delta x)\hat{\tau}_{P_{31}}(H_0(\delta x)-T)dx
    \end{eqnarray}
    equals \[\int_{G(F)\backslash G(\mathbb{A})^1}\sum_{n_1\in F^\times}\int_{\mathbb{A}^4}f(x^{-1}zu(0,n_1,n_2,n_3,n_4,n_5)x)\hat{\tau}_{P_{31}}(H_0(x)-T)dx.\]
\end{lemma}
\begin{proof}
It is clear that $M_{31,u(0,1,0,0,0,0)}\subset P_0\cap M_{31}$.

Write (\ref{+211}) as
\begin{eqnarray*}
    \int_{P_{211}(F)\backslash P_{211}(\mathbb{A})^1}\sum_{\delta\in M_{31,u_z}(F)\backslash P_{211}(F)}f_K(x^{-1}\delta^{-1}u_z\delta x)dx,
\end{eqnarray*}which is 
\begin{eqnarray}\label{sssd}
    \int_{P_{211}(F)\backslash P_{211}(\mathbb{A})^1}\sum_{(n_1,n_3)\in V^0}f_K(x^{-1}zu(0,n_1,0,n_3,0,0)x)dx.
\end{eqnarray}
By (\ref{mean1}), (\ref{sssd}) equals 
    \begin{eqnarray}\label{oop}
        \int_{P_{211}(F)\backslash P_{211}(\mathbb{A})}\int_{\mathbb{A}^2}f_{K}(p^{-1}zu(0,n_1,0,n_3,0,0)p)dn_1\,dn_3\,dp.
    \end{eqnarray}
    In fact, we have the equality
    \begin{eqnarray*}
        \rm{lim}_{T\rightarrow\infty}\int_{\mathbb{A}^2}f_{K}(p^{-1}zu(0,n_1,0,n_3,0,0)p)\hat{\tau}_{P_{22}}(H_0(p)-T)dn_1\,dn_3\,dp=(\ref{oop}).
    \end{eqnarray*}
    Write (\ref{oop}) as 
    \[\int_{M_{211}(F)\backslash M_{211}(\mathbb{A})}\int_{\mathbb{A}^5}\int_{\mathbb{A}^2}f_K(m^{-1}n^{-1}zu(0,n_1,0,n_3,0,0)nm)dn_1\,dn_3\,dn\,dm.\]
    Now we view $\rm{diag}(1,a,b,1)u(0,0,0,x,0,0)$ as the upper triangular matrix $P_{\rm{GL}(2)}$ of $\rm{GL}(2)$ and view $u(0,0,0,x,0,0)$ as the unipotent matrix $N_{\rm{GL}(2)}$ of $P_{\rm{GL}(2)}$. Write $Z_{\rm{GL}(2)}$ is the center of $\rm{GL}(2)$, and $K_{\rm{GL}(2)}$ as the maximal compact group of $\rm{GL}(2)$.

    Since the group $\{u(0,0,0,x,0,0),x\in\mathbb{A}\}$ is abelian, by the Poisson summation formula 
    \[\sum_{v\in N(\mu),v\neq e}\Phi_\mu(f:pvp^{-1})=\sum_{\xi\in X(\mu)_F,\xi\neq e}\Psi_\mu(\xi,p)+\Psi_\mu(0,p)-\Phi_\mu(f:e)\] in \cite[ p.~376]{A1}, where $p\in P_{\rm{GL}(2)}$, $\mu\in Z_{\rm{GL}(2)}(F)$, $N(\mu)=N_{\rm{GL}(2)}$, $X(\mu)$ is the dual space of $N_{\rm{GL}(2)}$, $e$ is the identity of $\rm{GL}(2)$, $\Phi_\mu(f:pvp^{-1})$ equals our $f_{K}$, and $\Psi_\mu(\xi,p)$ is the Fourier transform associated to $\{u(0,0,0,x,0,0),x\in\mathbb{A}\}$ (of course all these functions are viewed as the function on $\rm{GL}(2)$).

    We can conclude that (\ref{oop}) equals 
    \begin{eqnarray*}
        \int_{P_{211}(F)\backslash P_{211}(\mathbb{A})}\int_\mathbb{A}\sum_{n_3\in F^\times,n_3\neq e}f_{K}(p^{-1}zu(0,n_1,0,n_3,0,0)p)dn_1\,dp.
    \end{eqnarray*}This completes the proof, if we replace $x$ by $w_{12}x$ in (\ref{+211}). Since $w_{12}\in M_{211}\subset M_{31}$, \[\hat{\tau}_{P_{31}}(H_0(w_{12}x)-T)=\hat{\tau}_{P_{31}}(H_0(x)-T).\]
\end{proof}
Since $K_{G,\mathfrak{o}_{211,z}}(x,x)=\sum_{\delta\in G_{u_z}(F)\backslash G(F)}f(x^{-1}\delta^{-1} u_z\delta x)$, \begin{eqnarray*}
    G_{u_z}(F)=\left(\begin{matrix}
        a&n_0&n_1&n_2\\
        &b&n_3&n_4\\
        &&a&n_1\\
        &&&a
    \end{matrix}\right)\subset P_{211}(F),\quad n_0,n_1,n_2,n_3,n_4\in F,\ a,b\in F^\times.
\end{eqnarray*}
Consider 
\begin{eqnarray*}\label{V^0}
    \sum_{\delta \in G_{u_z}(F)\backslash P_{211}(F)}f(x^{-1}\delta^{-1}u_z\delta x)=\sum_{\substack{n_1,n_5\in F^\times}}\sum_{n_2,n_3,n_4\in F}f(x^{-1}zu(0,n_1,n_2,n_3,n_4,n_5)x).
\end{eqnarray*}
By the equality\[\sum_{\{P_3:P_3\supset P_1\}}(-1)^{\rm{dim}\, A_1\backslash A_3}\tau_1^3(H)\hat{\tau}_3(H)=0,\]decompose $J_{\mathfrak{o}_{211,z}}^T(f)$ as the sum 
\begin{eqnarray}\label{2111}
    &\int_{G(F)\backslash G(\mathbb{A})^1}[\sum_{\substack{n_1,n_5\in F^\times}}\sum_{n_2,n_3,n_4\in F}f(x^{-1}zu(0,n_1,n_2,n_3,n_4,n_5)x)\\\notag
    &-\sum_{n_1\in F^\times}\int_{\mathbb{A}^4}f(x^{-1}zu(0,n_1,n_2,n_3,n_4,n_5)x)\tau^{P_{31}}_{P_{211}}(H_0(x)-T)dn_2dn_3dn_4dn_5\\\notag
   &-\sum_{n_5\in F^\times}\int_{\mathbb{A}^4}f(x^{-1}zu(0,n_1,n_2,n_3,n_4,n_5)x)\tau_{P_{211}}^{P_{22}}(H_0(x)-T)dn_1dn_2dn_3dn_4\\\notag
    &+\int_{N_{211}(\mathbb{A})}f(x^{-1}zu(0,n_1,n_2,n_3,n_4,n_5)x)dn\tau_{P_{211}}(H_0(x)-T)]dx
\end{eqnarray}
\begin{eqnarray}\label{2112}
    &+\int_{G(F)\backslash G(\mathbb{A})^1}[\sum_{n_5\in F^\times}\int_{\mathbb{A}^4}f(x^{-1}zu(0,n_1,n_2,n_3,n_4,n_5)x)dn_1dn_2dn_3dn_4\\\notag
    &-\int_{N_{211}(\mathbb{A})}f(x^{-1}zu(0,n_1,n_2,n_3,n_4,n_5)x)dn\tau_{P_{211}}^{P_{22}}(H_0(x)-T)\\\notag
    &(\tau_{P_{211}}^{P_{31}}(H_0(x)-T)-\hat{\tau}_{P_{22}}(H_0(x)-T))]dx
\end{eqnarray}
\begin{eqnarray}\label{2113}
    &+\int_{G(F)\backslash G(\mathbb{A})^1}[\sum_{n_1\in F^\times}\int_{\mathbb{A}^4}f(x^{-1}zu(0,n_1,n_2,n_3,n_4,n_5)x)dn_2dn_3dn_4dn_5\\\notag
    &-\int_{N_{211}(\mathbb{A})}f(x^{-1}zu(0,n_1,n_2,n_3,n_4,n_5)x)dn\tau_{P_{211}}^{P_{31}}(H_0(x)-T)\\\notag
    &(\tau_{P_{211}}^{P_{22}}(H_0(x)-T)-\hat{\tau}_{P_{31}}(H_0(x)-T))]dx.
\end{eqnarray}
Write the Langlands decomposition\[x=nm\rm{diag}(a^{-1},a^{-1},1,d)k\in G(F)\backslash G(\mathbb{A})^1,\]where \[n\in N_{211}(F)\backslash N_{211}(\mathbb{A}),\ m\in M_{211}(F)\backslash M_{211}(\mathbb{A}),\ a,d\in F^\times\backslash\mathbb{A}^\times,\ k\in K,\] the modular function \[\delta^{-1}_{P_{211}}(a)=\rm{exp}(-\langle 2\rho_{P_{211}},H_0(a)\rangle)=|a|^4|d|^3.\]

Put \[F_z(x,y)=\int_{\mathbb{A}^3}f_K(zu(0,x,n_2,n_3,n_4,y))dn_2dn_3dn_4,\]
Thus (\ref{2111}) becomes
\begin{eqnarray*}
    &\int_{(F^\times \backslash\mathbb{A}^\times)^2}\int_{N_{211}(F)\backslash N_{211}(\mathbb{A})}\Big[\sum_{n_1,n_5\in F^\times}|ad|F_z(an_1,dn_5)\\
    &-\chi(|d|<\rm{exp}(T_1+T_3-2T_2))\sum_{n_1\in F^\times}|a|F_z^2(an_1,0)\\
    &-\chi(|a|<\rm{exp}(T_2-2T_3))\sum_{n_5\in F^\times}|d|F_z^1(0,dn_5)\\
    &+\chi(|d|<\rm{exp}(T_3-2T_2),|a|<\rm{exp}(T_2-2T_3))F_z^{1,2}(0,0)\Big]d^\times ad^\times d,
\end{eqnarray*}
where $\chi$ denotes the characteristic function, and $F^i,i=1,2$ denotes the Fourier transform of $F$ in its $i$-th argument, with respect to the standard additive character.

Define $J_{\mathfrak{o}_{211,z}}^T(f,\lambda)$ as the product of $J_{\mathfrak{o}_{211,z}}^T(f)$ and $\rm{exp}(-\langle\lambda,H_0(x)\rangle)$, where $\lambda\in\mathbb{C}^2$.
Write $\lambda=\lambda_2\hat{\alpha}_2+\lambda_3\hat{\alpha}_3$, by the Langlands decomposition, we have 
\begin{eqnarray*}\label{2inS}
    \rm{exp}(-\langle\lambda,H_0(x)\rangle)=|a|^{\frac{2\lambda_2+\lambda_3}{2}}|d|^{\frac{2\lambda_2+3\lambda_3}{4}}.
\end{eqnarray*}
Set \[\tilde{f}_{K,z}(n_1,n_3,n_5)=\int_{\mathbb{A}^2}f_K(zu(0,0,n_2,0,n_4,0)u(n_1,0,n_3,0,n_5)dn_2dn_4,\]
and\[(\frac{2\lambda_2+\lambda_3}{2}+1,\frac{2\lambda_2+3\lambda_3}{4}+1)=(s_2,s_3),\quad (2T_3-T_2,-T_1-T_3+2T_2)=(T_2',T_3').\]
Recall that if provided both limits of $f(s)$ and $g(s)$ exist as $s$ approaches $1$, \[\rm{lim}_{s\rightarrow1}f(s)g(s)=\rm{lim}_{s\rightarrow1}f(s)\cdot\rm{lim}_{s\rightarrow1}g(s).\]
Since \[|a|^{\frac{2\lambda_2+\lambda_3}{2}}=m^{\frac{2\lambda_2+\lambda_3}{2}}\cdot||(a,a)||^{\frac{2\lambda_2+\lambda_3}{2}},\] where $m$ is a constant depends on if $F=\mathbb{Q}$, we can replace $|a|^{\frac{2\lambda_2+\lambda_3}{2}}$ by $||(a,a)||^{\frac{2\lambda_2+\lambda_3}{2}}$.

We need to consider 
\begin{eqnarray}\label{pole211}
    \rm{lim}_{(s_2,s_3)\rightarrow (1,1)}\Pi_{i=2}^3\frac{\partial}{\partial s_i}\Pi_{j=2}^{3}(s_j-1)J_{\mathfrak{o}_{211,z}}^T(f,\lambda).
\end{eqnarray}
(\ref{2111}) in (\ref{pole211}) equals
\begin{eqnarray}\label{zeta2111}
    \rm{lim}_{(s_2,s_3)\rightarrow (1,1)}\Pi_{i=2}^3\frac{\partial}{\partial s_i}\Pi_{j=2}^{3}(s_j-1)\frac{\rm{vol}_{M_{211}}}{2c_F^2}\zeta_A(\tilde{f}_{K,z},(s_1,s_2),(T_2',T_3')).
\end{eqnarray}
Denote $\mathcal{J}(i_1,i_2)$ to be the function of 
\[\int_{F_S^5}f_{K_S}(zu(0,n_1,n_2,n_3,n_4,n_5))(\rm{log}||(n_1,n_3)||_S)^{i_1}(\rm{log}|n_5|_S)^{i_2}dn_1dn_2dn_3dn_4dn_5,\]where $i_1,i_2\in \mathbb{Z}_{\ge 0}$.

If we take $T=T_0=(0,0,0)$, (\ref{zeta2111}) equals 
\begin{eqnarray*}
    \frac{\rm{vol}_{M_{211}}}{2}\mathcal{J}(1,1)+\frac{\rm{vol}_{M_{211}}}{2c_F}\mathfrak{c}_F(S)\{\mathcal{J}(0,1)+\mathcal{J}(1,0)\}+\frac{\rm{vol}_{M_{211}}}{2c_F^2}\mathfrak{c}_F(S)^2\mathcal{J}(0,0).
\end{eqnarray*}
Write \[|n_{5}|=|n_5|_S|n_5|_{\overline{S}},\]where $\overline{S}$ is the complement of $S$.

For each $v\notin S$, \[\int_{\mathcal{O}_v}\rm{log}|n_{5,v}|_vdn_{5,v}=\frac{q_v^{-2}\rm{log}(q_v^{-1})}{1-q_v^{-2}},\quad \int_{\mathcal{O}_v}(\rm{log}|n_{5,v}|_v)^2dn_{5,v}=\frac{(q_v^{2}+1)(\rm{log} q_v)^2}{(q_v^{2}-1)^2}\]
and\[\rm{vol}_{M_{211}}\int_{N_{211}(\mathbb{A})}f_{K_S}(zn)dn.\]
Also, we can easily obatin that 
\[\frac{\frac{d^2}{ds^2}\zeta_F^S(s)|_{s=2}}{\zeta_F^S(2)}=\frac{(q_v^{2}+1)(\rm{log}q_v)^2}{(q_v^{2}-1)^2}.\]
As for the term corresponding to (\ref{2112}), we need to consider
\[\int_{\mathbb{A}^\times}|a|^{\frac{2\lambda_2+\lambda_3}{2}}(\tau_{P_{211}}^{P_{31}}(H_0(x)-T)-\hat{\tau}_{P_{22}}(H_0(x)-T))d^\times a,\]which equals \[c_F(\frac{2\lambda_2+\lambda_3}{2})^{-1}(\rm{exp}((T_1+T_3-2T_2)\frac{2\lambda_2+\lambda_3}{2})-|d|^{-\frac{1}{2}\frac{2\lambda_2+\lambda_3}{2}}\rm{exp}(-T_2\frac{2\lambda_2+\lambda_3}{2})),\]if $2\lambda_2+\lambda_3\neq 0$, and \[c_F(\frac{\rm{log}|d|}{2}+T_1+T_2-T_3),\]if $2\lambda_2+\lambda_3=0.$

It is clear that the integral \[\int_{N_{211}(\mathbb{A})}\sum_{n_5\in F^\times}\int_{\mathbb{A}^4}f(n^{-1}zu(0,n_1,n_2,n_3,n_4,n_5)n)dn_1dn_2dn_3dn_4dn\]equals \[\int_{N_{22}(\mathbb{A})}\sum_{n_5\in F^\times}\int_{\mathbb{A}^4}f(n^{-1}zu(0,n_1,n_2,n_3,n_4,n_5)n)dn_1dn_2dn_3dn_4dn,\]
since the group $\{n(0,0,0,0,0,a),\ a\in F\}$ is abelian. Set \[\tilde{f}_{K,z,5}(n_5)=\int_{N_{22}(\mathbb{A})}f_K(znu(0,0,0,0,n_5))dn,\]
then (\ref{2112}) equals \[\frac{\rm{vol}_{M_{211}}}{2c_F}\frac{\rm{exp}(-T_3'(s_2-1))\zeta_{\{2\}}(\tilde{f}_{K,z,5},s_3,T_3')-\rm{exp}(-T_2(s_3-1))\zeta_{\{2\}}(\tilde{f}_{K,z,5},s_2-\frac{s_3-1}{2},T_3')}{s_2-1}.\]
Thus if we take $T=T_0=(0,0,0)$, (\ref{2112}) in (\ref{pole211}) equals 
\begin{eqnarray*}
    \frac{\rm{vol}_{M_{211}}}{2c_F}\rm{lim}_{(s_2,s_3)\rightarrow (1,1)}\frac{\partial}{\partial s_3}(s_3-1)\{\frac{\partial}{\partial s_2}\zeta_{\{2\}}(\tilde{f}_{K,z,5},s_2-\frac{s_3-1}{2},T_0)\},
\end{eqnarray*}
which is
\begin{eqnarray*}
    \frac{\rm{vol}_{M_{211}}}{8}\mathcal{J}(0,2)+\frac{\rm{vol}_{M_{211}}}{4}\frac{q_v^{-2}\rm{log}(q_v^{-1})}{1-q_v^{-2}}\mathcal{J}(0,1)+\frac{\rm{vol}_{M_{211}}}{4}\frac{(q_v^{2}+1)(\rm{log}q_v)^2}{(q_v^{2}-1)^2}\mathcal{J}(0,0).
\end{eqnarray*}
As for the term corresponding to (\ref{2113}), we need to consider \[\int_{\mathbb{A}^\times}|d|^{\frac{2\lambda_2+3\lambda_3}{4}}(\tau_{P_{211}}^{P_{22}}(H_0(x)-T)-\hat{\tau}_{P_{31}}(H_0(x)-T))d^\times d,\]which equals \[c_F(\frac{2\lambda_2+3\lambda_3}{4})^{-1}(\rm{exp}(T_2-2T_3)\frac{2\lambda_2+3\lambda_3}{4}-|a|^{-\frac{2}{3}\frac{2\lambda_2+3\lambda_3}{4}}\rm{exp}(-\frac{4T_2}{3}\frac{2\lambda_2+3\lambda_3}{4})),\]if $2\lambda_2+3\lambda_3\neq 0$, and \[c_F(\frac{2\rm{log}|a|}{3}+T_2-\frac{2T_3}{3}),\]if $2\lambda_2+3\lambda_3=0$.

Set\[\tilde{f}_{K,z,1}(n_1,n_3)=\int_{\mathbb{A}^3}f_K(zu(0,n_2,0,n_4,n_5)u(n_1,0,n_3,0,0))dn_2\,dn_4\,dn_5.\]
Since \[\langle-2\rho_{P_{211}},H_0(x)\rangle=\langle(-2,-2,1,3),(-\rm{log}|a|,-\rm{log}|a|,0,\rm{log}|d|)\rangle=4\rm{log}|a|+3\rm{log}|d|,\]
we replace $|a|$ by $|a|^{\frac{4}{3}}$ in (\ref{2113}), and we need to change corresponding $s_2$.

By the formula of $\zeta_A^T(\phi,s)$ we obtained in (\ref{zetaa}), (\ref{2113}) equals the product of $\frac{3}{4}$ and  
\begin{eqnarray*}
    \frac{\rm{vol}_{M_{211}}}{2c_F}\frac{\rm{exp}(-T_2'(s_3-1))\zeta_{\{1\}}(\tilde{f}_{K,z,1},s_2,T_2')-\rm{exp}(-\frac{4T_2}{3}(s_3-1))\zeta_{\{1\}}(\tilde{f}_{K,z,1},s_2-\frac{2}{3}(s_3-1),T_2')}{s_3-1}.
\end{eqnarray*}
Write \[||(n_{1},n_{3})||=||(n_{1},n_{3})||_S||(n_{1},n_{3})||_{\overline{S}},\]where $\overline{S}$ is the complement of $S$.

For each $v\notin S$, \[\int_{\mathcal{O}_v^2}\rm{log}||(n_{1,v},n_{3,v})||_vdn_{1,v}dn_{3,v}=\frac{q_v^{-2}\rm{log}(q_v^{-1})}{1-q_v^{-2}},\]
\[\int_{\mathcal{O}_v^2}(\rm{log}||(n_{1,v},n_{3,v})||_v)^2dn_{1,v}dn_{3,v}=\frac{(q_v^{2}+1)(\rm{log}q_v)^2}{(q_v^{2}-1)^2}\]
and\[\rm{vol}_{M_{211}}\int_{N_{211}(\mathbb{A})}f_{K_S}(zn)dn.\]
Also, we can easily obatin that 
\[\frac{\frac{d^2}{ds^2}\zeta_F^S(s)|_{s=2}}{\zeta_F^S(2)}=\frac{(q_v^{2}+1)(\rm{log}q_v)^2}{(q_v^{2}-1)^2}.\]
Thus if we take $T=T_0=(0,0,0)$, (\ref{2113}) in (\ref{pole211}) equals 
\begin{eqnarray*}
    \frac{3}{4}\cdot\frac{\rm{vol}_{M_{211}}}{2c_F}\rm{lim}_{(s_2,s_3)\rightarrow (1,1)}\frac{\partial}{\partial s_2}(s_2-1)\{\frac{\partial}{\partial s_3}\zeta_{\{1\}}(\tilde{f}_{K,z,1},s_2-\frac{2}{3}(s_3-1),T_0)\},
\end{eqnarray*}
which is
\begin{eqnarray*}
    \frac{\rm{vol}_{M_{211}}}{8}\mathcal{J}(2,0)+\frac{\rm{vol}_{M_{211}}}{4}\frac{\frac{d}{ds}\zeta_F^S(s)|_{s=2}}{\zeta_F^S(2)}\mathcal{J}(1,0)+\frac{\rm{vol}_{M_{211}}}{4}\frac{\frac{d^2}{ds^2}\zeta_F^S(s)|_{s=2}}{\zeta_F^S(2)}\mathcal{J}(0,0).
\end{eqnarray*}
Thus we have
\begin{lemma}\label{mmm}
    For any $f\in C_c^\infty(G(\mathbb{A})^1)$, 
    \begin{eqnarray*}
        \rm{lim}_{(s_2,s_3)\rightarrow (1,1)}\Pi_{i=2}^3\frac{\partial}{\partial s_i}\Pi_{j=2}^{3}(s_j-1)J_{\mathfrak{o}_{211,z}}^T(f,\lambda)
    \end{eqnarray*} is the sum 
    \begin{eqnarray*}
    &\frac{\rm{vol}_{M_{211}}}{2}\mathcal{J}(1,1)+\frac{\rm{vol}_{M_{211}}}{2c_F}\mathfrak{c}_F(S)\{\mathcal{J}(0,1)+\mathcal{J}(1,0)\}+\frac{\rm{vol}_{M_{211}}}{2c_F^2}\mathfrak{c}_F(S)^2\mathcal{J}(0,0)\\
    &+\frac{\rm{vol}_{M_{211}}}{8}\mathcal{J}(0,2)+\frac{\rm{vol}_{M_{211}}}{4c_F}\frac{\frac{d}{ds}\zeta_F^S(s)|_{s=2}}{\zeta_F^S(2)}\mathcal{J}(0,1)+\frac{\rm{vol}_{M_{211}}}{4c_F}\mathfrak{c}_F'(S)\mathcal{J}(0,0)\\
    &+\frac{\rm{vol}_{M_{211}}}{8}\mathcal{J}(2,0)+\frac{\rm{vol}_{M_{211}}}{4c_F}\frac{\frac{d}{ds}\zeta_F^S(s)|_{s=2}}{\zeta_F^S(2)}\mathcal{J}(1,0)+\frac{\rm{vol}_{M_{211}}}{4c_F}\mathfrak{c}_F'(S)\mathcal{J}(0,0).
\end{eqnarray*}
\end{lemma}
\subsection{$\mathfrak{o}_{1,z}$ in $J_{\mathfrak{o}_{1111,z}}^T(f)$}

Recall that $J_{\mathfrak{o}_{1,z}}^T(f)$ equals the integral over $G(F)\backslash G(\mathbb{A})^1$ of
\begin{eqnarray*}
    &K_{G,\mathfrak{o}_{1,z}}(x,x)+\sum_{P\neq G}(-1)^{\rm{dim}\,Z\backslash A}\sum_{\delta\in P(F)\backslash G(F)}K_{P,\{\rm{La}_M(zu(1,0,0,1,0,1))\}_M}(\delta x,\delta x)\\
    &\cdot\hat{\tau}_P(H_0(\delta x)-T).
\end{eqnarray*}
Set $u_z=zu(1,0,0,1,0,1)$, for $z\in Z(F)$.

Since $G_{u_z}\subset P_0$, the integral over $G(F)\backslash G(\mathbb{A})^1$ of $K_{G,\mathfrak{o}_{1,z}}(x,x)$ is \[\int_{P_0(F)\backslash G(\mathbb{A})^1}\sum_{\delta\in G_{u_z}(F)\backslash P_0(F)}f(x^{-1}\delta^{-1}u_z\delta x)dx,\]
which is 
\[\int_{P_0(F)\backslash G(\mathbb{A})^1}\sum_{n_1,n_4,n_6\in F^\times}\sum_{n_2,n_3,n_5\in F}f(x^{-1}zu(n_1,n_2,n_3,n_4,n_5,n_6)x)dx.\]
Put $S_1=\{1,4,6\}.$ For any subset $S'\subset S_1$, we can find a unique standard parabolic subgroup $P_{S'}$ whose simple roots are $\{\alpha_{\lceil\frac{i}{2}\rceil}: i\in S'\}$.

According to the equality \[\sum_{\{P_3:P_3\supset P_0\}}(-1)^{\rm{dim}\,A_0\backslash A_3}\tau_0^3(H)\hat{\tau}_3(H)=0,\]
we can write $J_{\mathfrak{o}_{1,z}}(f)$ as the sum of \begin{eqnarray}\label{o11}
    &\int_{P_0(F)\backslash G(\mathbb{A})^1}\{\sum_{n_1,n_4,n_6\in F^\times}\sum_{n_2,n_3,n_5\in F}f(x^{-1}u(n_1,n_2,n_3,n_4,n_5,n_6)x)\\\notag
    &-\sum_{S'\subsetneq S_1}(-1)^{|S'|}\sum_{n_i\in F^\times,i\in S'}\int_{\mathbb{A}^{6-|S'|}}f(x^{-1}u(n_1,n_2,n_3,n_4,n_5,n_6)x)dn\\\notag
    &\tau_{P_0}^{P_{S'}}(H_0(x)-T)\}dx
\end{eqnarray}
and 
\begin{eqnarray}\label{o12}
    &\sum_{S'\subsetneq S_1}\int_{P_0(F)\backslash G(\mathbb{A})^1}\{\sum_{S''\subset S'}(-1)^{|S'|-|S''|}\sum_{n_i\in F^\times,i\in S''}\\\notag
    &\int_{\mathbb{A}^{6-|S''|}}f(x^{-1}zu(n_1,n_2,n_3,n_4,n_5,n_6)x)dn\tau_{P_0}^{P_{S_1-S'+S''}}(H_0(x)-T)\}\\\notag
    &\cdot(\tau_{P_0}^{P_{S'}}(H_0(x)-T)-\hat{\tau}_{P_{S_1-S'}}(H_0(x)-T))dx.
\end{eqnarray}
Put 
\begin{eqnarray*}
    F_z(x_1,x_2,x_3)=\int_{\mathbb{A}^3}f_K(\left(\begin{matrix}
        z&x_1&n_2&n_3\\
        &z&x_2&n_5\\
        &&z&x_3\\
        &&&z
    \end{matrix}\right))dn_2dn_3dn_5.
\end{eqnarray*}
Write $\lambda=\lambda_1\hat{\alpha}_1+\lambda_2\hat{\alpha}_2+\lambda_3\hat{\alpha}_3$ and suppose we have the Iwasawa decomposition \[x=u(n_{12},n_{13},n_{14},n_{23},n_{24},n_{34})\rm{diag}(a_1^{-1}a_2^{-1},a_2^{-1},1,a_3)k\in P_0(F)\backslash G(\mathbb{A})^1,\]the modular function is $|a_1|^3|a_2|^4|a_3|^3$.

Write \[(\frac{3\lambda_1+2\lambda_2+\lambda_3}{4}+1,\frac{\lambda_1+2\lambda_2+\lambda_3}{2}+1,\frac{\lambda_1+2\lambda_2+3\lambda_3}{4}+1)=(s_1,s_2,s_3)\]and \[(2T_1-T_2,2T_3-T_1-T_3,2T_3-T_2)=(T_1',T_2',T_3').\]
Then (\ref{o11}) becomes 
\begin{eqnarray*}
    &\int_{(F^\times\backslash \mathbb{A})^3}\{\sum_{n_1,n_4,n_6\in F^\times}|a_1a_2a_3|F_z(a_1n_1,a_2n_4,a_3n_6)\\
    &-\sum_{S'\subsetneq S_1}(-1)^{|S'|}(\Pi_{i\in S'}\chi(a_{\lceil\frac{i}{2}\rceil}<-T'_{\lceil\frac{i}{2}\rceil}))\sum_{n_i\in F^\times,i\in S'}\Pi_{i\in S'}\\
    &|a_{\lceil\frac{i}{2}\rceil}|\Pi_{j\in S_1-S'}|a_{\lceil\frac{i}{2}\rceil}|^{-1}F_z^{\{S_1-S'\}}(a_1n_1,a_2n_4,a_3n_6)\}d^\times a_1d^\times a_2d^\times a_3,
\end{eqnarray*}
where $\chi$ denotes the characteristic function, and $F_z^{\{S_1-S'\}}(a_1n_1,a_2n_4,a_3n_6)$ denotes the constant term of the Fourier transform of $(\Pi_{i\in S_1-S'}|a_{\lceil\frac{i}{2}\rceil}|^{-1})F_z(a_1n_1,a_2n_4,a_3n_6)$ with respect to the $\lceil\frac{i}{2}\rceil$-position, for $i\in S_1-S'$.

Write $J^T_{\mathfrak{o}_{1,z}}(f,\lambda)$ as the product of the factor $\rm{exp}(-\langle\lambda,H_0(x)\rangle)$ and the sum of (\ref{o11}) and (\ref{o12}). By the Iwasawa decomposition, we have \[\rm{exp}(-\langle\lambda,H_0(x)\rangle)=|a_1|^{s_1-1}|a_2|^{s_2-1}|a_3|^{s_3-1}.\]
Set \[\tilde{f}_{K,z}(n)=\int_{\mathbb{A}^3}f_K(zu(0,n_2,n_3,0,n_5,0))dn_2dn_3dn_4.\]
And we need to consider every single term in $J_{\mathfrak{o}_{1,z}}^T(f,\lambda)$ after cancelling the poles.

The constant term of the Laurent expansion of (\ref{o11}) at $(s_1,s_2,s_3)=(1,1,1)$ equals 
\begin{eqnarray}\label{limo1}
    \rm{lim}_{(s_1,s_2,s_3)\rightarrow (1,1,1)}\Pi_{i=1}^3\frac{\partial}{\partial s_i}\Pi_{j=1}^{3}(s_j-1)\frac{\rm{vol}_{M_0}}{4c_F^3}\zeta_A(\tilde{f}_{K,z},(s_1,s_2,s_3),(T'_1,T'_2,T'_3)).
\end{eqnarray}
Denote $\mathcal{J}(i_1,i_2,i_3)$ to be the function of 
\[\int_{F_S^6}f_{K_S}(zu(n_1,n_2,n_3,n_4,n_5,n_6))\Pi_{j\in S_1}(\rm{log}|n_{j}|_S)^{i_{\lceil\frac{j}{2}\rceil}}dn.\]
If we take $T=T_0=(0,0,0,0)$, (\ref{limo1}) equals
\begin{eqnarray}\label{sum1}
    &\frac{\rm{vol}_{M_0}}{4}\mathcal{J}(1,1,1)+\frac{\rm{vol}_{M_0}}{4c_F}\mathfrak{c}_F(S)\{\mathcal{J}(1,1,0)+\mathcal{J}(1,0,1)+\mathcal{J}(0,1,1)\}\\\notag
    &+\frac{\rm{vol}_{M_0}}{4c_F^2}\mathfrak{c}_F^2(S)\{\mathcal{J}(1,0,0)+\mathcal{J}(0,1,0)+\mathcal{J}(0,0,1)\}+\frac{\rm{vol}_{M_0}}{4c_F^3}\mathfrak{c}_F^3(S)\mathcal{J}(0,0,0).
\end{eqnarray}
Next, we decompose (\ref{o12}) as the sum 
\begin{eqnarray*}
    &\int_{P_0(F)\backslash G(\mathbb{A})}\{\sum_{n_1,n_4\in F^\times}\int_{\mathbb{A}^4}f-\sum_{n_1\in F^\times}\int_{\mathbb{A}^5}f\tau_{P_0}^{P_{121}}-\sum_{n_4\in F^\times}\int_{\mathbb{A}^5}f\tau_{P_0}^{P_{211}}+\int_{N_0(\mathbb{A})}f\tau_{P_0}^{P_{31}}\}\\
    &\cdot(\tau_{P_0}^{P_{112}}-\hat{\tau}_{P_{31}})dx\\
    &+\int_{P_0(F)\backslash G(\mathbb{A})}\{\sum_{n_1,n_6\in F^\times}\int_{\mathbb{A}^4}f-\sum_{n_1\in F^\times}\int_{\mathbb{A}^5}f\tau_{P_0}^{P_{112}}-\sum_{n_6\in F^\times}\int_{\mathbb{A}^5}f\tau_{P_0}^{P_{211}}+\int_{N_0(\mathbb{A})}f\tau_{P_0}^{P_{22}}\}\\
    &\cdot(\tau_{P_0}^{P_{121}}-\hat{\tau}_{P_{22}})dx\\
    &+\int_{P_0(F)\backslash G(\mathbb{A})}\{\sum_{n_4,n_6\in F^\times}\int_{\mathbb{A}^4}f-\sum_{n_4\in F^\times}\int_{\mathbb{A}^5}f\tau_{P_0}^{P_{112}}-\sum_{n_6\in F^\times}\int_{\mathbb{A}^5}f\tau_{P_0}^{P_{121}}+\int_{N_0(\mathbb{A})}f\tau_{P_0}^{P_{13}}\}\\
    &\cdot(\tau_{P_0}^{P_{211}}-\hat{\tau}_{P_{13}})dx\\
    &-\int_{P_0(F)\backslash G(\mathbb{A})}\{\sum_{n_1\in F^\times}\int_{\mathbb{A}^5}f-\int_{N_0(\mathbb{A})}f\tau_{P_0}^{P_{211}}\}(\tau_{P_0}^{P_{13}}-\hat{\tau}_{P_{211}})\\
    &-\int_{P_0(F)\backslash G(\mathbb{A})}\{\sum_{n_4\in F^\times}\int_{\mathbb{A}^5}f-\int_{N_0(\mathbb{A})}f\tau_{P_0}^{P_{121}}\}(\tau_{P_0}^{P_{22}}-\hat{\tau}_{P_{121}})\\
    &-\int_{P_0(F)\backslash G(\mathbb{A})}\{\sum_{n_6\in F^\times}\int_{\mathbb{A}^5}f-\int_{N_0(\mathbb{A})}f\tau_{P_0}^{P_{112}}\}(\tau_{P_0}^{P_{31}}-\hat{\tau}_{P_{112}}),
\end{eqnarray*}
where $f$ represents $f(x^{-1}zu(n_1,n_2,n_3,n_4,n_5,n_6)x)$, $\tau_{P_0}^P$ and $\hat{\tau}_P$ represent $\tau_{P_0}^P(H_0(x)-T)$ and $\hat{\tau}_P(H_0(x)-T)$, for any standard parabolic $P$.

For the first term, we need to consider \[\int_{\mathbb{A}^\times}|a_3|^{s_3-1}(\tau_{P_0}^{P_{112}}-\hat{\tau}_{P_{31}})d^\times a_3.\]
After cancelling the pole,  the result of the first term is 
\begin{eqnarray}\label{sum2}
   \notag&\frac{\rm{vol}_{M_0}}{24}\mathcal{J}(2,1,0)+\frac{\rm{vol}_{M_0}}{12c_F}\mathfrak{c}_F(S)\mathcal{J}(1,1,0)+\frac{\rm{vol}_{M_0}}{24c_F}\mathfrak{c}_F(S)\mathcal{J}(2,0,0)+\frac{\rm{vol}_{M_0}}{12c^2_F}\mathfrak{c}^2_F(S)\mathcal{J}(1,0,0)\\\notag
   &+\frac{\rm{vol}_{M_0}}{12c_F^2}\mathfrak{c}'_F(S)c_F^S\mathcal{J}(0,1,0)+\frac{\rm{vol}_{M_0}}{12c^3_F}\mathfrak{c}'_F(S)\mathfrak{c}_F(S)c_F^S\mathcal{J}(0,0,0)\\\notag
   &+\frac{\rm{vol}_{M_0}}{12}\mathcal{J}(1,2,0)+\frac{\rm{vol}_{M_0}}{6c_F}\mathfrak{c}_F(S)\mathcal{J}(1,1,0)+\frac{\rm{vol}_{M_0}}{12c_F}\mathfrak{c}_F(S)\mathcal{J}(0,2,0)+\frac{\rm{vol}_{M_0}}{6c^2_F}\mathfrak{c}^2_F(S)\mathcal{J}(0,1,0)\\
   &+\frac{\rm{vol}_{M_0}}{6c_F^2}\mathfrak{c}'_F(S)c_S\mathcal{J}(1,0,0)+\frac{\rm{vol}_{M_0}}{6c^3_F}\mathfrak{c}'_F(S)\mathfrak{c}_F(S)c_S\mathcal{J}(0,0,0).
\end{eqnarray}
The second and third terms are similar, which are 
\begin{eqnarray}\label{sum3}
    \notag&\frac{\rm{vol}_{M_0}}{16}\mathcal{J}(2,0,1)+\frac{\rm{vol}_{M_0}}{8c_F}\mathfrak{c}_F(S)c_S\mathcal{J}(1,0,1)+\frac{\rm{vol}_{M_0}}{16c_F}\mathfrak{c}_F(S)\mathcal{J}(2,0,0)+\frac{\rm{vol}_{M_0}}{8c^2_F}\mathfrak{c}^2_F(S)\mathcal{J}(1,0,0)\\\notag
    &+\frac{\rm{vol}_{M_0}}{8c_F^2}\mathfrak{c}'_F(S)c_F^S\mathcal{J}(0,0,1)+\frac{\rm{vol}_{M_0}}{8c^3_F}\mathfrak{c}'_F(S)\mathfrak{c}_F(S)c_F^S\mathcal{J}(0,0,0)
    \\\notag
    &+\frac{\rm{vol}_{M_0}}{16}\mathcal{J}(1,0,2)+\frac{\rm{vol}_{M_0}}{8c_F}\mathfrak{c}_F(S)\mathcal{J}(1,0,1)+\frac{\rm{vol}_{M_0}}{16c_F}\mathfrak{c}_F(S)\mathcal{J}(0,0,2)+\frac{\rm{vol}_{M_0}}{8c^2_F}\mathfrak{c}^2_F(S)\mathcal{J}(0,0,1)\\
    &+\frac{\rm{vol}_{M_0}}{8c_F^2}\mathfrak{c}'_F(S)c_F^S\mathcal{J}(1,0,0)+\frac{\rm{vol}_{M_0}}{8c^3_F}\mathfrak{c}'_F(S)\mathfrak{c}_F(S)c_F^S\mathcal{J}(0,0,0),
\end{eqnarray}
and 
\begin{eqnarray}\label{sum4}
    &\frac{\rm{vol}_{M_0}}{12}\mathcal{J}(0,2,1)+\frac{\rm{vol}_{M_0}}{6c_F}\mathfrak{c}_F(S)c_S\mathcal{J}(0,1,1)+\frac{\rm{vol}_{M_0}}{12c_F}\mathfrak{c}_F(S)c_S\mathcal{J}(0,2,0)\\\notag
    &+\frac{\rm{vol}_{M_0}}{6c^2_F}\mathfrak{c}^2_F(S)c_S\mathcal{J}(0,1,0)+\frac{\rm{vol}_{M_0}}{6c_F^2}\mathfrak{c}'_F(S)c_S\mathcal{J}(0,0,1)+\frac{\rm{vol}_{M_0}}{6c^3_F}\mathfrak{c}'_F(S)\mathfrak{c}_F(S)c_S\mathcal{J}(0,0,0)
    \\\notag
    &+\frac{\rm{vol}_{M_0}}{24}\mathcal{J}(0,1,2)+\frac{\rm{vol}_{M_0}}{12c_F}\mathfrak{c}_F(S)c_S\mathcal{J}(0,1,1)+\frac{\rm{vol}_{M_0}}{24c_F}\mathfrak{c}_F(S)c_S\mathcal{J}(0,0,2)\\\notag
    &+\frac{\rm{vol}_{M_0}}{12c^2_F}\mathfrak{c}^2_F(S)c_S\mathcal{J}(0,0,1)+\frac{\rm{vol}_{M_0}}{12c_F^2}\mathfrak{c}'_F(S)c_S\mathcal{J}(0,1,0)+\frac{\rm{vol}_{M_0}}{12c^3_F}\mathfrak{c}'_F(S)\mathfrak{c}_F(S)c_S\mathcal{J}(0,0,0).
\end{eqnarray}
The term \[\int_{(\mathbb{A}^\times)^2}|a_2|^{s_2-1}|a_3|^{s_3-1}\tau_{P_0}^{P_{13}}d^\times a_2d^\times a_3\] equals zero, if we take $T=T_0$.

The term \[\int_{(\mathbb{A}^\times)^2}|a_2|^{s_2-1}|a_3|^{s_3-1}\hat{\tau}_{P_{211}}d^\times a_2d^\times a_3\]equals 
\[2\frac{|a_1|^{-\frac{1}{2}s_2+\frac{3}{4}}}{(\frac{3}{2}s_2-s_3-1)(-s_2+2s_3-1)}.\]
To cancel the poles, take $s_3=1$, multiply by $(s_2-1)^2$,  then take the derivative with respect to $s_2$ twice. The result of the fourth term is 
\begin{eqnarray}\label{sum5}
    \notag&\frac{\rm{vol}_{M_0}}{144}\mathcal{J}(3,0,0)+\frac{\rm{vol}_{M_0}}{48c_F}\mathfrak{c}_F(S)\mathcal{J}(2,0,0)\\
    &+\frac{\rm{vol}_{M_0}}{24c_F^2}\mathfrak{c}_F'(S)c_F^S\mathcal{J}(1,0,0)+\frac{\rm{vol}_{M_0}}{24c_F^3}\mathfrak{c}_F''(S)(c_F^S)^2\mathcal{J}(0,0,0),
\end{eqnarray}
Similar is the fifth and the sixth, which are 
\begin{eqnarray}\label{sum6}
    \notag&\frac{\rm{vol}_{M_0}}{36}\mathcal{J}(0,3,0)+\frac{\rm{vol}_{M_0}}{12c_F}\mathfrak{c}_F(S)\mathcal{J}(0,2,0)\\
    &+\frac{\rm{vol}_{M_0}}{6c^2_F}\mathfrak{c}_F'(S)c_F^S\mathcal{J}(0,1,0)+\frac{\rm{vol}_{M_0}}{6c^3_F}\mathfrak{c}_F''(S)(c_F^S)^2\mathcal{J}(0,0,0),
\end{eqnarray}and
\begin{eqnarray}\label{sum7}
    \notag&\frac{\rm{vol}_{M_0}}{144}\mathcal{J}(0,0,3)+\frac{\rm{vol}_{M_0}}{48c_F}\mathfrak{c}_F(S)\mathcal{J}(0,0,2)\\
    &+\frac{\rm{vol}_{M_0}}{24c^2_F}\mathfrak{c}_F'(S)c_F^S\mathcal{J}(0,0,1)+\frac{\rm{vol}_{M_0}}{24c^3_F}\mathfrak{c}_F''(S)(c_F^S)^2\mathcal{J}(0,0,0).
\end{eqnarray}
Collecte all contributions according to the number of logarithmic factors, we obtain, we have 
\begin{lemma}\label{6.5}
    For $z\in Z(F)$ and $f\in C_c^\infty(G(\mathbb{A})^1)$, $J_{\mathfrak{o}_{1,z}}^{T_0}(f)$ equals the sum
    \begin{eqnarray*}
        (\ref{sum1})+(\ref{sum2})+(\ref{sum3})+(\ref{sum4})+(\ref{sum5})+(\ref{sum6})+(\ref{sum7}).
    \end{eqnarray*}
\end{lemma}

\section{The coefficient $a^M(S,u)$ }\label{8}
The goal of this section is to compute explicitly the coefficients $a^M(S,u)$ appearing in the geometric expansion of the trace formula (\ref{ogamma}).
\begin{lemma}\label{v_P}
    For the element \begin{eqnarray*}
        n=\left(\begin{matrix}
            1&n_1&n_2&n_3\\
            &1&n_4&n_5\\
            &&1&n_6\\
            &&&1
        \end{matrix}\right)\in N_0(\mathbb{A}),
    \end{eqnarray*}then $v_{sP_0}(\lambda,n,T_0)$ is given by Table \ref{table1}.
    
\end{lemma}
\begin{proof}
    Recall 
    \begin{eqnarray*}
        v_{sP_0}(\lambda,n,T_0)=e^{\langle\lambda,s^{-1}H_0(w_sn)\rangle}.
    \end{eqnarray*}
    We illustrate the computation by taking $s=(13)(24)$, we have
    \begin{eqnarray*}
        w_sn=\left(\begin{matrix}
            &&1&n_6\\
            &&&1\\
            1&n_1&n_2&n_3\\
            &1&n_3&n_5
        \end{matrix}\right).
    \end{eqnarray*}Then by the definition of $||\cdot||$,
    \begin{eqnarray*}
        &a_4=||(1,n_4,n_5)||,\quad a_3a_4=||(1,n_4,n_5,n_2-n_1n_4,n_3-n_1n_5,n_2n_5-n_3n_4)||,\\
        &a_2a_3a_4=||(1,n_1)||,\quad a_1a_2a_3a_4=1,
    \end{eqnarray*}we can write 
    \begin{eqnarray*}
        &s^{-1}H_0(w_sn)=(\rm{log}\frac{||(1,n_4,n_5,n_2-n_1n_4,n_3-n_1n_5,n_2n_5-n_3n_4)||}{||(1,n_4,n_5)||},\rm{log}||(1,n_4,n_5)||,\\
        &-\rm{log}||(1,n_4,n_2-n_1n_4)||,\rm{log}\frac{||(1,n_4,n_2-n_1n_4)||}{||(1,n_4,n_5,n_2-n_1n_4,n_3-n_1n_5,n_2n_5-n_3n_4)||}).
    \end{eqnarray*}
    Since $\lambda=\lambda_1\hat{\alpha}_1+\lambda_2\hat{\alpha}_2+\lambda_3\hat{\alpha}_3$, the factor \begin{eqnarray*}
        v_{sP_0}(\lambda,n,T_0)=e^{\langle\lambda,s^{-1}H_0(w_sn)\rangle}.
    \end{eqnarray*} equals 
    \begin{eqnarray*}
        &||(1,n_4,n_5,n_2-n_1n_4,n_3-n_1n_5,n_2n_5-n_3n_4)||^{-\lambda_1-\lambda_2-\lambda_3}||(1,n_4,n_5)||^{\lambda_1}\\
        &||(1,n_4,n_2-n_1n_4)||^{\lambda_3}.
    \end{eqnarray*}
    The remaining cases can be obtained by analogous computations. 
\end{proof}
\begin{lemma}\label{7.2}
    For \begin{eqnarray*}
        n_0=\left(\begin{matrix}
            1&v_1&v_2&v_3\\
            &1&v_4&v_5\\
            &&1&v_6\\
            &&&1
        \end{matrix}\right)\in N_0(F_S),
    \end{eqnarray*}then
    \begin{eqnarray*}
        &w_{M_0}(1,n,T_0)=\frac{1}{6}\{(\rm{log}\,|v_1|)^3+4(\rm{log}\,|v_4|)^3+(\rm{log}\,|v_6|)^3+6(\rm{log}\,|v_1|)^2\rm{log}\,|v_4|\\
        &+9(\rm{log}\,|v_1|)^2\rm{log}\,|v_6|+12\rm{log}\,|v_1|(\rm{log}\,|v_4|)^2+12(\rm{log}\,|v_4|)^2\rm{log}\,|v_6|+9\rm{log}\,|v_1|(\rm{log}\,|v_6|)^2\\
        &+6\rm{log}\,|v_4|(\rm{log}\,|v_6|)^2+36\rm{log}\,|v_1|\rm{log}\,|v_4|\rm{log}\,|v_6|\}.
    \end{eqnarray*}
\end{lemma}
\begin{proof}
    For this $n_0$, 
    \begin{eqnarray*}
        n_0^{-1}=\left(\begin{matrix}
        1&-v_1&v_1v_4-v_2&v_1v_5+v_2v_6-v_3-v_1v_4v_6\\
        &1&-v_4&v_4v_6-v_5\\
        &&1&-v_6\\
        &&&1
    \end{matrix}\right),
    \end{eqnarray*}and for $a_0=\rm{diag}(a.b,c,1)$, set $v'$ equals
    \begin{eqnarray*}
        a_0^{-1}n_0^{-1}a_0n_0=\left(\begin{matrix}
            1&\frac{a-b}{a}v_1&\frac{a-c}{a}v_2+\frac{c-b}{a}v_1v_4&\frac{a-1}{a}v_3+\frac{1-b}{a}v_1v_5+\frac{c-1}{a}v_1v_4v_6+\frac{1-c}{a}v_2v_6\\
            &1&\frac{b-c}{b}v_4&\frac{b-1}{b}v_5+\frac{1-c}{b}v_4v_6\\
            &&1&\frac{c-1}{c}v_6\\
            &&&1
        \end{matrix}\right).
    \end{eqnarray*}
    Set\begin{eqnarray*}
        v'=\left(\begin{matrix}
            1&\xi_1&\xi_2&\xi_3\\
            &1&\xi_4&\xi_5\\
            &&1&\xi_6\\
            &&&1
        \end{matrix}\right),
    \end{eqnarray*}we have
    \begin{eqnarray*}
        v_1=\frac{a}{a-b}\xi_1,\quad v_2=\frac{a}{a-c}\xi_2+\frac{ab}{(a-b)(a-c)}\xi_1\xi_4,
    \end{eqnarray*}
    \begin{eqnarray*}
        v_3=\frac{a}{a-1}\xi_3+\frac{ab}{(a-1)(a-b)}\xi_1\xi_5+\frac{ac}{(a-1)(a-c)}\xi_2\xi_6+\frac{abc}{(a-1)(a-b)(a-c)}\xi_1\xi_4\xi_6,
    \end{eqnarray*}
    \begin{eqnarray*}
        v_4=\frac{b}{b-c}\xi_4,\quad v_5=\frac{b}{b-1}\xi_5+\frac{bc}{(b-1)(b-c)}\xi_4\xi_6,\quad v_6=\frac{c}{c-1}\xi_6.
    \end{eqnarray*}
    Recall that 
    \[r_\beta(\lambda,1,a_0)=|(a_0^\beta-1)(1-a_0^{-\beta})|^{-\frac{\lambda(\beta_1)}{2}},\]where $\beta_1\in\mathfrak{a}_0$ is a coroot scaled by a nonnegative factor in such a way that the limits $w_{P_0}(\lambda,1,v)$ exist and
    are nonzero for generic $\lambda$ and $v$.
    \[w_{sP_0}(\lambda,1,v,T_0)=v_{P_0}(\lambda,n,T_0)\cdot\Pi_{\beta}r_\beta(\lambda,1,a_0),\]$\beta$ runs through the reduced roots of $(N_{P_0} /(N_{P_0}\cap N_{w_s^{-1}P_0w_s}), A_{M_0})$.

    We take $s=(13)(24)$ for example. In this case, $\beta\in\{\alpha_2,\alpha_1+\alpha_2,\alpha_2+\alpha_3,\alpha_1+\alpha_2+\alpha_3\}$. Then $\Pi_\beta r_\beta(\lambda,1,a_0)$ equals
    \begin{eqnarray*}
        (\frac{(b-c)^2}{|bc|})^{-\frac{\lambda((\alpha_2)_1)}{2}}(\frac{(a-c)^2}{|ac|})^{-\frac{\lambda((\alpha_1+\alpha_2)_1)}{2}}(\frac{(b-1)^2}{|b|})^{-\frac{\lambda((\alpha_2+\alpha_3)_1)}{2}} (\frac{(a-1)^2}{|a|})^{-\frac{\lambda((\alpha_1+\alpha_2+\alpha_3)_1)}{2}}.
    \end{eqnarray*}
    By the replacement of $v_i$ by $\xi_i$, $i\in\{1,2,3,4,5,6\}$ in $v_{sP_0}(\lambda,v,T_0)$, we can write $v_{sP_0}(\lambda,v,T_0)$ as the product of 
    \begin{eqnarray*}
        &||(1,\frac{b}{b-c}\xi_4,\frac{b}{b-1}\xi_5+\frac{bc}{(b-1)(b-c)}\xi_4\xi_6,\frac{a}{a-c}\xi_2+\frac{ab}{(a-c)(b-c)}\xi_1\xi_4,\frac{a}{a-1}\xi_3\\
        &+\frac{ab}{(a-1)(b-1)}\xi_1\xi_5+\frac{ac}{(a-1)(a-c)}\xi_2\xi_6,\frac{ab^2c}{(a-1)(b-1)(a-c)(b-c)}\xi_1\xi_4^2\xi_6)||^{-\lambda_1-\lambda_2-\lambda_3}
    \end{eqnarray*}and
    \begin{eqnarray*}
        ||(1,\frac{b}{b-c}\xi_4,\frac{b}{b-1}\xi_5+\frac{bc}{(b-1)(b-c)}\xi_4\xi_6)||^{\lambda_1}||(1,\frac{b}{b-c}\xi_4,\frac{a}{a-c}\xi_2+\frac{ab}{(a-c)(b-c)}\xi_1\xi_4)||^{\lambda_3}.
    \end{eqnarray*}
    To make sure the limit $w_{sP_0}(\lambda,1,v)$ exists, it follows that \[(\alpha_2)_1=\alpha_2,(\alpha_1+\alpha_2)_1=\alpha_1+\alpha_2,(\alpha_2+\alpha_3)_1=\alpha_2+\alpha_3,(\alpha_1+\alpha_2+\alpha_3)_1=\alpha_1+\alpha_2+\alpha_3.\]
    Then by Lemma \ref{v_P} and similar calculations, we have the every term in Table \ref{table2}.

    Since $w_{sP_0}(\lambda,1,v,T_0)$ is a $(G,M_0)$-family, by Lemma 3.4 in \cite{A2}, we can take $\lambda_1=\lambda_2=-1$, and $\lambda_3=-1.5$, then \[\lambda=-\hat{\alpha}_1-\hat{\alpha}_2-1.5\hat{\alpha}_3=(-\frac{13}{8},-\frac{5}{8},\frac{3}{8},\frac{15}{8}).\]By Lemma 3.4 in \cite{A2}, since \[w_{M_0}(1,n,T_0)=\sum_{s\in S_4}\frac{e^{\langle\lambda,w_{sP_0}(\lambda,1,v,T_0)\rangle}}{\Pi_{\eta\in\Phi_P}\langle\lambda,s^{-1}\eta\rangle},\]
    it follows that 
    \begin{eqnarray*}
        &w_{M_0}(1,n,T_0)=\frac{1}{6}\{(\rm{log}\,|v_1|)^3+4(\rm{log}\,|v_4|)^3+(\rm{log}\,|v_6|)^3+6(\rm{log}\,|v_1|)^2\rm{log}\,|v_4|\\
        &+9(\rm{log}\,|v_1|)^2\rm{log}\,|v_6|+12\rm{log}\,|v_1|(\rm{log}\,|v_4|)^2+12(\rm{log}\,|v_4|)^2\rm{log}\,|v_6|+9\rm{log}\,|v_1|(\rm{log}\,|v_6|)^2\\
        &+6\rm{log}\,|v_4|(\rm{log}\,|v_6|)^2+36\rm{log}\,|v_1|\rm{log}\,|v_4|\rm{log}\,|v_6|\}.
    \end{eqnarray*}(For the complete calculation, use the code in Appendix \ref{code}.)
\end{proof}
\begin{lemma}\label{211v}
    For the element 
    \begin{eqnarray*}
        n_{211}=\left(\begin{matrix}
            1&&n_1&n_2\\
            &1&n_3&n_4\\
            &&1&n_5\\
            &&&1
        \end{matrix}\right)\in N_{211}(\mathbb{A}),
    \end{eqnarray*}then $v_{sP_{211}}(\lambda,n_{211},T_0)$ equals
    \begin{eqnarray*}
        \begin{cases}
            1,&\text{if}\ s=1,\\
            ||(1,n_5)||^{-\lambda_3},&\text{if}\ s=(34),\\
            ||(1,n_1,n_3)||^{-\lambda_2},&\text{if}\ s=(13),\\
            ||(1,n_5)||^{\lambda_2}||(1,n_5,n_4-n_3n_5,n_2-n_1n_5)||^{-\lambda_2-\lambda_3},&\text{if}\ s=(143),\\
            ||(1,n_1,n_3)||^{\lambda_3}||(1,n_1,n_2,n_3,n_4,n_1n_4-n_2n_3)||^{-\lambda_2-\lambda_3},&\text{if}\ s=(13)(24),\\
            ||(1,n_1,n_2,n_3,n_4,n_1n_4-n_2n_3)||^{-\lambda_2}||(1,n_5,n_4-n_3n_5,n_2-n_1n_5)||^{-\lambda_3},&\text{if}\ s=(14)(23),
        \end{cases}
    \end{eqnarray*}where $\lambda=\lambda_2\hat{\alpha}_2+\lambda_3\hat{\alpha}_3$.
\end{lemma}
\begin{lemma}
    For 
    \begin{eqnarray*}
        v_{211}=\left(\begin{matrix}
            1&&v_1&v_2\\
            &1&v_3&v_4\\
            &&1&v_5\\
            &&&1
        \end{matrix}\right)\in N_{211}(F_S),\quad v_0=\left(\begin{matrix}
            1&v_1&v_2&v_3\\
            &1&v_4&v_5\\
            &&1&v_6\\
            &&&1
        \end{matrix}\right)\in N_0(F_S),
    \end{eqnarray*}then $w_{M_{211}}(1,v_{211},T_0)$ equals
    \begin{eqnarray*}
        \frac{1}{2}\{(\rm{log}||(v_1,v_3)||)^2+4\rm{log}|v_5|\rm{log}||(v_1,v_3)||+(\rm{log}|v_5|)^2\},
    \end{eqnarray*}and $w_{M_{211}}(1,v_{0},T_0)$ equals 
    \begin{eqnarray*}
        \frac{1}{2}\{(\rm{log}|v_1|)^2+2(\rm{log}|v_4|)^2+(\rm{log}|v_6|)^2+4\rm{log}|v_1|\rm{log}|v_4|+4\rm{log}|v_4|\rm{log}|v_6|+4\rm{log}|v_1|\rm{log}|v_6|\}.
    \end{eqnarray*}
\end{lemma}
\begin{proof}
    For this $v_{211}$,
    \begin{eqnarray*}
        v_{211}^{-1}=\left(\begin{matrix}
            1&&-v_1&v_1v_5-v_2\\
            &1&-v_3&v_3v_5-v_4\\
            &&1&-v_5\\
            &&&1
        \end{matrix}\right),
    \end{eqnarray*}and for $a_0=\rm{diag}(a,a,1,d)$,
    \begin{eqnarray*}
        v'=a_0^{-1}v_{211}^{-1}a_0v_{211}=\left(\begin{matrix}
            1&&\frac{a-1}{a}v_1&\frac{a-d}{a}v_2+\frac{d-1}{a}v_1v_5\\
            &1&\frac{a-1}{a}v_3&\frac{a-d}{a}v_4+\frac{d-1}{a}v_3v_5\\
            &&1&(1-d)v_5\\
            &&&1
        \end{matrix}\right).
    \end{eqnarray*}
    Set \begin{eqnarray*}
        v'=\left(\begin{matrix}
            1&&\xi_1&\xi_2\\
            &1&\xi_3&\xi_4\\
            &&1&\xi_5\\
            &&&1
        \end{matrix}\right),
    \end{eqnarray*}we have
    \begin{eqnarray*}
        v_1=\frac{a}{a-1}\xi_1,\quad v_2=\frac{a}{a-d}\xi_2+\frac{a}{(a-1)(a-d)}\xi_1\xi_5,\quad v_3=\frac{a}{a-1}\xi_3,\\
        v_4=\frac{a}{a-d}\xi_4+\frac{a}{(a-d)(a-1)}\xi_3\xi_5,\quad v_5=\frac{1}{1-d}\xi_5.
    \end{eqnarray*}Then by Lemma \ref{211v},
    \begin{eqnarray*}
        w_{sP_{211}}(\lambda,1,v_{211},T_0)=\begin{cases}
            1,&\text{if}\ s=1,\\
            |v_5|_S^{-\lambda_3},&\text{if}\ s=(34),\\
            ||(v_1,v_3)||_S^{-\lambda_2},&\text{if}\ s=(13),\\
            ||(v_1,v_3)||_S^{-\lambda_2-\lambda_3}|v_5|_S^{-\lambda_3},&\text{if}\ s=(143),\\
            ||(v_1v_5,v_3v_5)||_S^{-\lambda_2-\lambda_3}||(v_1,v_3)||_S^{\lambda_3},&\text{if}\ s=(13)(24),\\
            ||((v_1v_5,v_3v_5)||_S^{-\lambda_2}||(v_3v_5,v_1v_5)||_S^{-\lambda_3},&\text{if}\ s=(14)(23),
        \end{cases}
    \end{eqnarray*}and
    \begin{eqnarray*}
        w_{sP_{211}}(\lambda,1,v_0,T_0)=\begin{cases}
            1,&\text{if}\ s=1,\\
            |v_6|_S^{-\lambda_3},&\text{if}\ s=(34),\\
            |v_1v_4|_S^{-\lambda_2},&\text{if}\ s=(13),\\
            |v_1v_4|_S^{-\lambda_2}|v_1v_4v_6|^{-\lambda_3},&\text{if}\ s=(143),\\
            |v_1v_4^2v_6|_S^{-\lambda_2}|v_4v_6|_S^{-\lambda_3},&\text{if}\ s=(13)(24),\\
            |v_1v_4^2v_6|_S^{-\lambda_2}|v_1v_4v_6|_S^{-\lambda_3},&\text{if}\ s=(14)(23).
        \end{cases}
    \end{eqnarray*}
    Take $\lambda=-\hat{\alpha}_2-1.5\hat{\alpha}_3$, by Lemma 3.4 in \cite{A2}, $w_{M_{211}}(1,v_{211},T_0)$ equals 
    \begin{eqnarray*}
        \frac{1}{2}\{\frac{(\frac{3}{2}\rm{log}|v_5|)^2}{-\frac{15}{4}}+\frac{(\rm{log}||(v_1,v_3)||)^2}{-\frac{5}{2}}+\frac{(\frac{5}{2}\rm{log}||(v_1,v_3)||+\frac{3}{2}\rm{log}|v_5|)^2}{-\frac{5}{2}}\\+\frac{\rm{log}||(v_1,v_3)||+\frac{5}{2}\rm{log}|v_5|)^2}{-\frac{15}{4}}+\frac{(\frac{5}{2}\rm{log}||(v_1,v_3||+\frac{5}{2}\rm{log}|v_5|)^2}{\frac{3}{2}}\},
    \end{eqnarray*}which is 
    \begin{eqnarray*}
        \frac{1}{2}\{(\rm{log}||(v_1,v_3)||)^2+4\rm{log}|v_5|\rm{log}||(v_1,v_3)||+(\rm{log}|v_5|)^2\}.
    \end{eqnarray*}And $w_{M_{211}}(1,v_{0},T_0)$ equals
    \begin{eqnarray*}
        &\frac{1}{2}\{\frac{(\frac{3}{2}\rm{log}|v_6|)^2}{-\frac{15}{4}}+\frac{(\rm{log}|v_1|+\rm{log}|v_4|)^2}{-\frac{5}{2}}+\frac{(\frac{5}{2}\rm{log}|v_1|+\frac{5}{2}\rm{log}|v_4|+\frac{3}{2}\rm{log}|v_6|)^2}{-\frac{5}{2}}\\
        &\frac{(\rm{log}|v_1|+\frac{7}{2}\rm{log}|v_4|+\frac{5}{2}\rm{log}|v_6|)^2}{-\frac{15}{4}}+\frac{(\frac{5}{2}\rm{log}|v_1|+\frac{7}{2}\rm{log}|v_4|+\frac{5}{2}\rm{log}|v_6|)^2}{\frac{3}{2}}\},
    \end{eqnarray*}which is 
    \begin{eqnarray*}
        \frac{1}{2}\{(\rm{log}|v_1|)^2+2(\rm{log}|v_4|)^2+(\rm{log}|v_6|)^2+4\rm{log}|v_1|\rm{log}|v_4|+4\rm{log}|v_4|\rm{log}|v_6|+4\rm{log}|v_1|\rm{log}|v_6|\}.
    \end{eqnarray*}
\end{proof}
\begin{theorem}\label{theorem1}
    For $\mathfrak{o}=\mathfrak{o}_{31,z}$, $z\in Z(F)$ and $f\in C_c^\infty(G(\mathbb{A})^1)$, $J_{\mathfrak{o}_{31,z}}(f)$ equals
    \[\rm{vol}_{M_{31}}J_{M_{31}}(1,f)+\rm{vol}_{M_{31}}\cdot \frac{\frac{d}{ds}\zeta_F^S(s)|_{s=3}}{\zeta_F^S(3)}J_G(u(0,0,1,0,0,0),f),\]that is \[a^{M_{31}}(S,1)=\rm{vol}_{M_{31}},\quad a^{G}(S,u(0,0,1,0,0,0))=\rm{vol}_{M_{31}}\cdot \frac{\frac{d}{ds}\zeta_F^S(s)|_{s=3}}{\zeta_F^S(3)}.\]
\end{theorem}
\begin{proof}
    Recall that $J_{\mathfrak{o}_{31,z}}(f)$ equals the sum of 
    \begin{eqnarray*}
        \rm{vol}_{M_{31}}\int_{F_S^3}f_{K_S}(zu(0,0,n_{14},0,n_{24},n_{34}))(\rm{log}||(n_{14},n_{24},n_{34})||)dn_{14}dn_{24}dn_{34}
    \end{eqnarray*}and \[\rm{vol}_{M_{31}}\cdot \frac{\frac{d}{ds}\zeta_F^S(s)|_{s=3}}{\zeta_F^S(3)}\cdot\int_{N_{31}(\mathbb{A})}f_{K_S}(zn)dn.\] 
    By calculation, we have 
    \[J_{M_{31}}(z,f)=\rm{vol}_{M_{31}}\int_{F_S^3}f_{K_S}(zu(0,0,n_{14},0,n_{24},n_{34}))(\rm{log}||(n_{14},n_{24},n_{34})||)dn_{14}dn_{24}dn_{34},\]and 
    \[J_G(u(0,0,1,0,0,0),f)=\int_{N_{31}(\mathbb{A})}f_{K_S}(zn)dn.\]
    This theorem follows.
\end{proof}

\begin{theorem}\label{theorem2}
    For $\mathfrak{o}=\mathfrak{o}_{22,z}$, $z\in Z(F)$ and $f\in C_c^\infty(G(\mathbb{A})^1)$, $J_{\mathfrak{o}_{22,z}}(f)$ equals
    \[\rm{vol}_{M_{22}}J_{M_{22}}(1,f)+\rm{vol}_{M_{22}}\cdot \mathfrak{c}_{2,F}(S)J_G(u(0,1,0,0,1,0),f),\]that is \[a^{M_{22}}(S,1)=\rm{vol}_{M_{22}},\quad a^{G}(S,u(0,1,0,0,1,0))=\rm{vol}_{M_{22}}\cdot \mathfrak{c}_{2,F}(S).\]
\end{theorem}
\begin{proof}
    Recall that $J_{\mathfrak{o}_{22,z}}(f)$ equals 
    \begin{eqnarray*}
    &\rm{vol}_{\rm{GL}(2)}\int_{\rm{GL}(2,F_S)}f_{K_S}(zu^{\rm{GL}(2)}(M_S))|\rm{det}(M_S)|^2\rm{log}|\rm{det}M_S|dM_S\cdot c^S_{2,F}\\
    &+\rm{vol}_{\rm{GL}(2)}\mathfrak{c}_{2,F}(S)\int_{\rm{GL}(2,F_S)}f_{K_S}(zu^{\rm{GL}(2)}(M_S))|\rm{det}(M_S)|^{2}dM_S.
\end{eqnarray*}
We have
\begin{eqnarray*}
    J_{M_{22}}(z,f)=\int_{N_{22}(F_S)}f_{K_S}(zn)w_{M_{22}}(1,n,T_0)dn,
\end{eqnarray*}
and 
\begin{eqnarray*}
    J_{G}(zu(0,1,0,0,1,0),f)=\int_{N_{22}(F_S)}f_{K_S}(zn)dn.
\end{eqnarray*}
By the fact that \[\rm{vol}_{M_{22}}=\rm{vol}_{\rm{GL}(2)}^2,\quad c_{2,S}d^\times g_S=|\rm{det}\,g_S|^2_Sd^+g_S,\quad \rm{vol}_{M_{22}}=c_{2,F}=c_{2,S}c_{2,F}^S,\]the theorem follows.
\end{proof}
Set 
    \begin{eqnarray*}
        u_{22}=\left(\begin{matrix}
            1&&&\\
            &1&&\\
            &&1&v_5\\
            &&&1
        \end{matrix}\right),\quad v_{22}=\left(\begin{matrix}
            1&&v_1&v_2\\
            &1&v_3&v_4\\
            &&1&\\
            &&&1
        \end{matrix}\right),
    \end{eqnarray*}and 
    \begin{eqnarray*}
        u_{31}=\left(\begin{matrix}
            1&&v_1&\\
            &1&&\\
            &&1&\\
            &&&1
        \end{matrix}\right),\quad v_{31}=\left(\begin{matrix}
            1&&&v_2\\
            &1&&v_4\\
            &&1&v_5\\
            &&&1
        \end{matrix}\right).
    \end{eqnarray*}
\begin{theorem}\label{theorem3}
    For $\mathfrak{o}=\mathfrak{o}_{211,z}$, $z\in Z(F)$ and $f\in C_c^\infty(G(\mathbb{A})^1)$, $J_{\mathfrak{o}_{211,z}}(f)$ equals
    \begin{eqnarray*}
        &\rm{vol}_{M_{211}}J_{M_{211}}(1,f)+\rm{vol}_{M_{211}}\frac{\mathfrak{c}_F(S)}{2c_F}J_{M_{22}}(u(0,0,0,0,0,1),f)\\
        &+\rm{vol}_{M_{211}}\frac{\frac{d}{ds}\zeta_F^S(s)|_{s=2}}{\zeta_F^S(2)}J_{M_{31}}(u(0,1,0,0,0,0),f)+\rm{vol}_{M_{211}}\frac{\frac{d^2}{ds^2}\zeta_F^S(s)|_{s=2}}{\zeta_F^S(2)}J_G(u(0,1,0,0,0,1),f),
    \end{eqnarray*}that is 
    \begin{eqnarray*}
        &a^{M_{211}}(S,1)=\rm{vol}_{M_{211}},\quad a^{M_{22}}(S,u(0,0,0,0,0,1))=\rm{vol}_{M_{211}}\frac{\mathfrak{c}_F(S)}{2c_F},\\
        &a^{M_{31}}(S,u(0,1,0,0,0,0))=\rm{vol}_{M_{211}}\frac{\frac{d}{ds}\zeta_F^S(s)|_{s=2}}{\zeta_F^S(2)},\\
        &a^G(S,u(0,1,0,0,0,1))=\rm{vol}_{M_{211}}\frac{\frac{d^2}{ds^2}\zeta_F^S(s)|_{s=2}}{\zeta_F^S(2)}.
    \end{eqnarray*}
\end{theorem}
\begin{proof}
    By calculation, we have 
    \begin{eqnarray*}
        w_{M_{211}}(1,u_{22}v_{22},T_0)=\rm{log}|v_5|+\rm{log}||(v_1,v_3)||,
    \end{eqnarray*}and 
    \begin{eqnarray*}
        w_{M_{31}}(1,u_{31}v_{31},T_0)=\rm{log}|v_5|+\rm{log}||(v_1,v_3)||.
    \end{eqnarray*}
    \begin{eqnarray*}
        w_{M_{211}}(1,v_{22},T_0)=\rm{log}|v_2v_3-v_1v_4|.
    \end{eqnarray*}
    Then we can write 
    \begin{eqnarray*}
        J_{M_{22}}(zu(0,0,0,0,0,1),f)=\int_{N_{211}(F_S)}f_{K_S}(zu_{22}v_{22})w_{M_{22}}(1,u_{22}v_{22},T_0)dv,
    \end{eqnarray*}
    \begin{eqnarray*}
        J_{M_{31}}(zu(0,1,0,0,0,0),f)=\int_{N_{211}(F_S)}f_{K_S}(zu_{31}v_{31})w_{M_{31}}(1,u_{31}v_{31},T_0)dv.
    \end{eqnarray*}
    \begin{eqnarray*}
        J_G(zu(0,1,0,0,0,1),f)=\int_{N_{211}(F_S)}f_{K_S}(zv)dv.
    \end{eqnarray*}
    and 
    \begin{eqnarray*}
        J_{M_{22}}(zu(1,0,0,0,1),f)=\int_{N_0(F_S)}f_{K_S}(zv_{22})w_{M_{22}}(1,v_{22},T_0)dv.
    \end{eqnarray*}
    Thus by Lemma \ref{mmm}, we can write $J_{\mathfrak{o}_{211},z}^{T_0}(f)$ as the sum 
    \begin{eqnarray*}
        &\rm{vol}_{M_{211}}J_{M_{211}}(1,f)+\rm{vol}_{M_{211}}\frac{\mathfrak{c}_F(S)}{2c_F}J_{M_{22}}(u(0,0,0,0,0,1),f)\\
        &+\rm{vol}_{M_{211}}\frac{\frac{d}{ds}\zeta_F^S(s)|_{s=2}}{\zeta_F^S(2)}J_{M_{31}}(u(0,1,0,0,0,0),f)+\rm{vol}_{M_{211}}\frac{\frac{d^2}{ds^2}\zeta_F^S(s)|_{s=2}}{\zeta_F^S(2)}J_G(u(0,1,0,0,0,1),f),
    \end{eqnarray*}that is,
    \begin{eqnarray*}
        &a^{M_{211}}(S,1)=\rm{vol}_{M_{211}},\quad a^{M_{22}}(S,u(0,0,0,0,0,1))=\frac{\mathfrak{c}_F(S)}{c_F},\\
        &a^{M_{31}}(S,u(0,1,0,0,0,0))=\rm{vol}_{M_{211}}\frac{\frac{d}{ds}\zeta_F^S(s)|_{s=2}}{\zeta_F^S(2)},\\
        &a^G(S,u(0,1,0,0,0,1))=\rm{vol}_{M_{211}}\frac{\frac{d^2}{ds^2}\zeta_F^S(s)|_{s=2}}{\zeta_F^S(2)}.
    \end{eqnarray*}
\end{proof}
\begin{theorem}\label{theorem4}
    For $\mathfrak{o}=\mathfrak{o}_{1,z}$, $z\in Z(F)$ and $f\in C_c^\infty(G(\mathbb{A})^1)$, $J_{\mathfrak{o}_{1,z}}(f)$ equals
    \begin{eqnarray*}
        &\rm{vol}_{M_0}J_{M_0}(z,f)+\rm{vol}_{M_0}\frac{\mathfrak{c}_F(S)}{2c_F}J_{M_{211}}(zu(1,0,0,0,0,0),f)\\
        &+\frac{\rm{vol}_{M_0}}{3c_F^2}(\mathfrak{c}_F(S)^2+\mathfrak{c}'_F(S)c_F^S)J_{M_{31}}(zu(1,0,0,1,0,0),f)\\
        &+\frac{\rm{vol}_{M_0}}{4c_F^2}\mathfrak{c}_F^2(S)J_{M_{22}}(zu(1,0,0,0,0,1),f)\\
        &+\frac{\rm{vol}_{M_0}}{4c_F^3}(\mathfrak{c}_F^3(S)+3\mathfrak{c}'_F(S)\mathfrak{c}_F(S)c_F^S+\mathfrak{c}''_F(S)c_F^S)J_{M_0}(zu(1,0,0,1,0,1),f),
    \end{eqnarray*}that is 
    \begin{eqnarray*}
        &a^{M_{0}}(S,1)=\rm{vol}_{M_{0}},\quad a^{M_{211}}(S,u(1,0,0,0,0,0))=\rm{vol}_{M_0}\frac{\mathfrak{c}_F(S)}{2c_F},\\
        &a^{M_{31}}(S,u(1,0,0,1,0,0))=\frac{\rm{vol}_{M_0}}{3c_F^2}(\mathfrak{c}_F(S)^2+\mathfrak{c}'_F(S)c_F^S),\\
        &a^{M_{22}}(S,u(1,0,0,0,0,1))=\frac{\rm{vol}_{M_0}}{4c_F^2}\mathfrak{c}_F^2(S),\\
        &a^{G}(S,u(0,1,0,0,1,0))=\frac{\rm{vol}_{M_0}}{4c_F^3}(\mathfrak{c}_F^3(S)+3\mathfrak{c}'_F(S)\mathfrak{c}_F(S)c_F^S+\mathfrak{c}''_F(S)(c_F^S)^2).
    \end{eqnarray*}
\end{theorem}
\begin{proof}
    Set 
    \begin{eqnarray*}
        u_{211}=\left(\begin{matrix}
           1&v_1&&\\
           &1&&\\
           &&1&\\
           &&&1 
        \end{matrix}\right),\quad u_{31}=\left(\begin{matrix}
           1&v_1&&\\
           &1&v_4&\\
           &&1&\\
           &&&1 
        \end{matrix}\right),\quad u_{22}=\left(\begin{matrix}
           1&v_1&&\\
           &1&&\\
           &&1&v_6\\
           &&&1 
        \end{matrix}\right),
    \end{eqnarray*}and
    \begin{eqnarray*}
        v_{211}=\left(\begin{matrix}
           1&&v_2&v_3\\
           &1&v_4&v_5\\
           &&1&v_6\\
           &&&1 
        \end{matrix}\right),\quad v_{31}=\left(\begin{matrix}
           1&&&v_3\\
           &1&&v_5\\
           &&1&v_6\\
           &&&1 
        \end{matrix}\right),\quad v_{22}=\left(\begin{matrix}
           1&&v_2&v_3\\
           &1&v_4&v_5\\
           &&1&\\
           &&&1 
        \end{matrix}\right),
    \end{eqnarray*}
    Recall that 
    \begin{eqnarray*}
        &w_{M_0}(1,v_0,T_0)=\frac{1}{6}\{(\rm{log}\,|v_1|)^3+4(\rm{log}\,|v_4|)^3+(\rm{log}\,|v_6|)^3+6(\rm{log}\,|v_1|)^2\rm{log}\,|v_4|\\
        &+9(\rm{log}\,|v_1|)^2\rm{log}\,|v_6|+12\rm{log}\,|v_1|(\rm{log}\,|v_4|)^2+12(\rm{log}\,|v_4|)^2\rm{log}\,|v_6|+9\rm{log}\,|v_1|(\rm{log}\,|v_6|)^2\\
        &+6\rm{log}\,|v_4|(\rm{log}\,|v_6|)^2+36\rm{log}\,|v_1|\rm{log}\,|v_4|\rm{log}\,|v_6|\},
    \end{eqnarray*}and $w_{M_{211}}(u_{211}v_{211},T_0)$ equals 
    \begin{eqnarray*}
        \frac{1}{2}\{(\rm{log}|v_1|)^2+2(\rm{log}|v_4|)^2+(\rm{log}|v_6|)^2+4\rm{log}|v_1|\rm{log}|v_4|+4\rm{log}|v_4|\rm{log}|v_6|+4\rm{log}|v_1|\rm{log}|v_6|\}.
    \end{eqnarray*}
    By calculation, we have 
    \begin{eqnarray*}
        &w_{M_{31}}(1,u_{31}v_{31},T_0)=\rm{log}|v_1|+\rm{log}|v_4|+\rm{log}|v_6|,\\
        &w_{M_{22}}(1,u_{22}v_{22},T_0)=\rm{log}|v_1|+2\rm{log}|v_4|+\rm{log}|v_6|,
    \end{eqnarray*}
    Thus, by Lemma \ref{6.5}, we have 
    \begin{eqnarray*}
        &J_{\mathfrak{o}_{1,z}}(f)=\rm{vol}_{M_0}J_{M_0}(z,f)+\rm{vol}_{M_0}\frac{\mathfrak{c}_F(S)}{2c_F}J_{M_{211}}(zu(1,0,0,0,0,0),f)\\
        &+\frac{\rm{vol}_{M_0}}{3c_F^2}(\mathfrak{c}_F(S)^2+\mathfrak{c}'_F(S)c_F^S)J_{M_{31}}(zu(1,0,0,1,0,0),f)\\
        &+\frac{\rm{vol}_{M_0}}{4c_F^2}\mathfrak{c}_F^2(S)J_{M_{22}}(zu(1,0,0,0,0,1),f)\\
        &+\frac{\rm{vol}_{M_0}}{4c_F^3}(\mathfrak{c}_F^3(S)+3\mathfrak{c}'_F(S)\mathfrak{c}_F(S)c_F^S+\mathfrak{c}''_F(S)c_F^S)J_{M_0}(zu(1,0,0,1,0,1),f).
    \end{eqnarray*}
    We can then obtain the desired coefficient $a^{M_0}(S,u)$.
\end{proof}
\subsection{The coefficient $a^{M}(S,\gamma)$}\label{uniddd}
In this section, we give all the coefficient $a^{M}(S,u)$ in 
\begin{eqnarray*}
    J_\mathfrak{o}^T(f)=\sum_{M\in\mathscr{L}(M_0)}|W_{0}^M||W_0^G|^{-1}\sum_{u\in (\mathcal{U}_M(F))_{M,S}}a^M(S,u)J_M^T(u,f).
\end{eqnarray*}
For any $\gamma\in G(F)$, write the Jordan decomposition $\gamma=\sigma u$, where $\sigma$ represents the semisimple part and $u$ represents the unipotent part. $S$ is the finite set of $\Sigma$ which contains $\Sigma_\infty$.

For any Levi subgroup $M$, set 
\begin{eqnarray}
    \epsilon^M(\sigma)=\begin{cases}
        1,\quad &if\  A_{M_{\sigma}}=A_M,\\
        0,\quad &else.
    \end{cases}
\end{eqnarray}
Arthur has proved the descent formula for coefficients
\begin{eqnarray}\label{descent}
    a^M(S,\gamma)=\epsilon^M(\sigma)|\iota^M(\sigma)|^{-1}\sum_{\{u|\sigma u\sim \gamma\}}a^{M_{\sigma}}(S,u),
\end{eqnarray}
in Theorem 8.1 of \cite{A5}, where the sum is taken over all $u\in (\mathcal{U}_{M_\sigma}(F))_{M_{\sigma,S}}$ such that $\sigma u$ is $(M,S)$-equivalent to $\gamma$. $\iota^M(\sigma)$ is  $|M_{\sigma}^0(F)\backslash M_{\sigma}(F)|$, that is, the number of the connected components. In the case of $\rm{GL}(n)$, it equals $1$.

The formula (\ref{descent}) means for another ramified orbit $\mathfrak{o}\neq \mathfrak{o}_{1111}^4$, we can descent to the unipotent case, though some $a^M(S,\gamma)$ may equals $0$.

Now by our previous calculation, we can obtain all the coefficient $a^M(S,\gamma)$.

\newpage

\section{appendix}
\begin{appendix}

\section{The explicit expressions of \(v_{sP_0}(\lambda,n,T_0)\)}\label{table11}
\renewcommand{\arraystretch}{1.6} 
\begin{longtable}{l p{12cm}}
\toprule
\textbf{Permutation $s$} & \textbf{Expression for $v_{sP_0}(\lambda, n, T_0)$}\label{table1} \\
\midrule
\endfirsthead

\toprule
\textbf{Permutation $s$} & \textbf{Expression for $v_{sP_0}(\lambda, n, T_0)$ (Continued)} \\
\midrule
\endhead

\bottomrule
\endfoot

\bottomrule
\endlastfoot

$s = 1$ & $1$ \\

$s = (12)$ & $\|(1, n_3)\|^{-\lambda_1}$ \\

$s = (23)$ & $\|(1, n_4)\|^{-\lambda_2}$ \\

$s = (34)$ & $\|(1, n_6)\|^{-\lambda_3}$ \\

$s = (13)$ & $\|(1, n_1, n_2)\|^{-\lambda_1} \|(1, n_4, n_2 - n_1n_4)\|^{-\lambda_2}$ \\

$s = (14)$ & 
$\begin{aligned}[t]
    & \|(1, n_1, n_2, n_3)\|^{-\lambda_1 - \lambda_2} \|(1, n_1, n_6, n_1n_6, n_3 - n_2n_6)\|^{\lambda_2}\\
    & \|(1, n_4, n_4n_6 - n_5, n_3 - n_1n_5 - n_2n_6 + n_1n_4n_6)\|^{-\lambda_2 - \lambda_3}
\end{aligned}$ \\

$s = (24)$ & $\|(1, n_4, n_5)\|^{-\lambda_2} \|(1, n_6, n_5 - n_4n_6)\|^{-\lambda_3}$ \\

$s = (123)$ & $\|(1, n_1, n_2)\|^{-\lambda_1 - \lambda_2} \|(1, n_2)\|^{2\lambda_2}$ \\

$s = (132)$ & $\|(1, n_3)\|^{\lambda_1} \|(1, n_4, n_2 - n_1n_4)\|^{-\lambda_1 - \lambda_2}$ \\

$s = (134)$ & 
$\begin{aligned}[t]
    & \|(1, n_1, n_2, n_3)\|^{-\lambda_1 - \lambda_2 - \lambda_3} \|(1, n_4, n_2 - n_1n_4)\|^{-\lambda_2} \|(1, n_1, n_2)\|^{\lambda_2 + \lambda_3}
\end{aligned}$ \\

$s = (143)$ & 
$\begin{aligned}[t]
    & \|(1, n_1, n_6, n_1n_6, n_3 - n_2n_6)\|^{-\lambda_1} \|(1, n_6)\|^{\lambda_1 + \lambda_2} \\
    & \|(1, n_6, n_5 - n_4n_6, n_3 - n_2n_6 - n_1n_5 + n_1n_4n_6)\|^{-\lambda_2 - \lambda_3}
\end{aligned}$ \\

$s = (234)$ & $\|(1, n_4)\|^{\lambda_2} \|(1, n_4, n_5)\|^{-\lambda_2 - \lambda_3}$ \\

$s = (243)$ & $\|(1, n_6)\|^{\lambda_3} \|(1, n_6, n_5 - n_4n_6)\|^{-\lambda_2 - \lambda_3}$ \\

$s = (124)$ & 
$\begin{aligned}[t]
    & \|(1, n_3)\|^{\lambda_2 + \lambda_3} \|(1, n_1, n_6, n_2n_6, n_3 - n_2n_6)\|^{-\lambda_3} \|(1, n_1, n_2, n_3)\|^{-\lambda_1 - \lambda_2}
\end{aligned}$ \\

$s = (142)$ & 
$\begin{aligned}[t]
    & \|(1, n_6, n_5 - n_4n_6, n_3 - n_1n_5 - n_2n_6 + n_1n_4n_6)\|^{-\lambda_2 - \lambda_3} \\
    &\|(1, n_1, n_2, n_3)\|^{-\lambda_1 - \lambda_2 - \lambda_3} \|(1, n_6, n_5 - n_4n_6)\|^{\lambda_1 + \lambda_2}
\end{aligned}$ \\

$s = (12)(34)$ & $\|(1, n_2)\|^{-\lambda_1} \|(1, n_6)\|^{-\lambda_3}$ \\

$s = (13)(24)$ & 
$\begin{aligned}[t]
    & \|(1, n_4, n_5, n_2 - n_1n_4, n_3 - n_1n_5, n_2n_5 - n_3n_4)\|^{-\lambda_1 - \lambda_2 - \lambda_3} \|(1, n_4, n_5)\|^{\lambda_1}\\
    & \|(1, n_4, n_2 - n_1n_4)\|^{\lambda_3}
\end{aligned}$ \\

$s = (14)(23)$ & 
$\begin{aligned}[t]
    & \|(1, n_1, n_2, n_3)\|^{-\lambda_1}\|(1, n_4, n_5, n_2 - n_1n_4, n_3 - n_1n_5, n_3n_4 - n_2n_5)\|^{-\lambda_2} \\
    &\|(1, n_6, n_5 - n_4n_6, n_3 - n_2n_6 - n_1n_5 + n_1n_4n_6)\|^{-\lambda_3}
\end{aligned}$ \\

$s = (1234)$ & $\|(1, n_1, n_2, n_3)\|^{-\lambda_1 - \lambda_2 - \lambda_3} \|(1, n_2)\|^{\lambda_2} \|(1, n_1, n_2)\|^{\lambda_3}$ \\

$s = (1243)$ & 
$\begin{aligned}[t]
    & \|(1, n_1, n_6, n_1n_6, n_3 - n_2n_6)\|^{-\lambda_1 - \lambda_2 - \lambda_3}  \|(1, n_1)\|^{\lambda_2 + \lambda_3} \|(1, n_6)\|^{\lambda_1 + \lambda_2}
\end{aligned}$ \\

$s = (1324)$ & 
$\begin{aligned}[t]
&\|(1, n_1, n_2, n_3)\|^{-\lambda_1} \|(1, n_4, n_2 - n_1n_4)\|^{\lambda_3} \\
&\|(1, n_4, n_5, n_2 - n_1n_4, n_3 - n_1n_5, n_3n_4 - n_2n_5)\|^{-\lambda_2 - \lambda_3}
\end{aligned}$ \\

$s = (1342)$ & 
$\begin{aligned}[t]
    &\|(1, n_4, n_2 - n_1n_4)\|^{-\lambda_1 - \lambda_2} \|(1, n_1)\|^{\lambda_1 + \lambda_2 + \lambda_3} \|(1, n_4, n_5)\|^{-\lambda_2 - \lambda_3}
\end{aligned}$ \\

$s = (1423)$ & 
$\begin{aligned}[t]
    & \|(1, n_4, n_5)\|^{\lambda_1}\|(1, n_4, n_5, n_2 - n_1n_4, n_3 - n_1n_5, n_2n_5 - n_3n_4)\|^{-\lambda_1 - \lambda_2} \\
    &\|(1, n_6, n_5 - n_4n_6, n_3 - n_1n_5 - n_2n_6 + n_1n_4n_6)\|^{-\lambda_3}
\end{aligned}$ \\

$s = (1432)$ & 
$\begin{aligned}[t]
    & \|(1, n_6)\|^{\lambda_3}\|(1, n_6, n_5 - n_4n_6, n_3 - n_1n_5 - n_2n_6 + n_1n_4n_6)\|^{-\lambda_1 - \lambda_2 - \lambda_3} \\
    & \|(1, n_6, n_5 - n_4n_6)\|^{\lambda_1}
\end{aligned}$ \\

\end{longtable}

\section{The explicit expressions of $w_{sP_0}(\lambda,1,v,T_0)$}\label{table22}
\renewcommand{\arraystretch}{1.6} 
\begin{longtable}{l p{12cm}}
\toprule
\textbf{Permutation $s$} & \textbf{Expression for $w_{sP_0}(\lambda,1,v,T_0)$}\label{table2} \\
\midrule
\endfirsthead

\toprule
a\textbf{Permutation $s$} & \textbf{Expression for $w_{sP_0}(\lambda,1,v,T_0)$ (Continued)} \\
\midrule
\endhead

\bottomrule
\endfoot

\bottomrule
\endlastfoot

$s = 1$ & $1$ \\

$s = (12)$ & $|v_1|_S^{-\lambda_1}$ \\

$s = (23)$ & $|v_4|_S^{-\lambda_2}$ \\

$s = (34)$ & $|v_6|_S^{-\lambda_3}$ \\

$s = (13)$ & $|v_1|_S^{-\lambda_1} |v_1 v_4|_S^{-\lambda_2}$ \\

$s = (14)$ & $|v_1 v_4 v_6|_S^{-\lambda_1} |v_1 v_4 v_6|_S^{-\lambda_2} |v_1 v_4 v_6|_S^{-\lambda_3}$ \\

$s = (24)$ & $|v_4 v_6|_S^{-\lambda_2} |v_4 v_6|_S^{-\lambda_3}$ \\

$s = (123)$ & $|v_1 v_4|_S^{-\lambda_1} |v_4|_S^{-\lambda_2}$ \\

$s = (132)$ & $|v_1|_S^{-\lambda_1} |v_1 v_4|_S^{-\lambda_2}$ \\

$s = (134)$ & $|v_1 v_4 v_6|_S^{-\lambda_1 - \lambda_2} |v_6|_S^{-\lambda_3}$ \\

$s = (143)$ & $|v_1 v_4|_S^{-\lambda_1} |v_1 v_4|_S^{-\lambda_2} |v_1 v_4 v_6|_S^{-\lambda_3}$ \\

$s = (234)$ & $|v_4 v_6|_S^{-\lambda_2} |v_6|_S^{-\lambda_3}$ \\

$s = (243)$ & $|v_4|_S^{-\lambda_2} |v_4 v_6|_S^{-\lambda_3}$ \\

$s = (124)$ & $|v_1 v_4 v_6|_S^{-\lambda_1} |v_4 v_6|_S^{-\lambda_2} |v_4 v_6|_S^{-\lambda_3}$ \\

$s = (142)$ & $|v_1|_S^{-\lambda_1} |v_1 v_4 v_6|_S^{-\lambda_2} |v_1 v_4 v_6|_S^{-\lambda_3}$ \\

$s = (12)(34)$ & $|v_1|_S^{-\lambda_1} |v_6|^{-\lambda_3}$ \\

$s = (13)(24)$ & $|v_1 v_4|_S^{-\lambda_1} |v_1 v_4^2 v_6|_S^{-\lambda_2} |v_4 v_6|_S^{-\lambda_3}$ \\

$s = (14)(23)$ & 
$\begin{aligned}[t]
    & |v_1 v_4 v_6|_S^{-\lambda_1} |v_1 v_4^2 v_6|_S^{-\lambda_2}  |v_1 v_4 v_6|_S^{-\lambda_3}
\end{aligned}$ \\

$s = (1234)$ & $|v_1 v_4 v_6|_S^{-\lambda_1} |v_4 v_6|_S^{-\lambda_2} |v_6|_S^{-\lambda_3}$ \\

$s = (1243)$ & $|v_1 v_4|_S^{-\lambda_1} |v_4|_S^{-\lambda_2} |v_4 v_6|_S^{-\lambda_3}$ \\

$s = (1324)$ & 
$\begin{aligned}[t]
    & |v_1 v_4 v_6|_S^{-\lambda_1} |v_1 v_4^2 v_6|_S^{-\lambda_2} |v_4 v_6|^{-\lambda_3}
\end{aligned}$ \\

$s = (1342)$ & $|v_1|_S^{-\lambda_1} |v_1 v_4 v_6|_S^{-\lambda_2} |v_6|_S^{-\lambda_3}$ \\

$s = (1423)$ & 
$\begin{aligned}[t]
    & |v_1 v_4|_S^{-\lambda_1} |v_1 v_4^2 v_6|_S^{-\lambda_2}  |v_1 v_4 v_6|^{-\lambda_3}
\end{aligned}$ \\

$s = (1432)$ & $|v_1|_S^{-\lambda_1} |v_1 v_4|_S^{-\lambda_2} |v_1 v_4 v_6|_S^{-\lambda_3}$ \\

\end{longtable}

\section{code}\label{code}

\begin{verbatim}
var('a, b, c')

lam = [-13/8, -5/8, 3/8, 15/8]

G = SymmetricGroup(4)

def theta_new(s):
    res = 1
    for i in range(1, 4):
        idx1 = s(i) - 1
        idx2 = s(i+1) - 1
        res *= (lam[idx1] - lam[idx2])
    return res

terms = [
    (a^3, "(1,2)"),
    (b^3, "(2,3)"),
    ((27/8)*c^3, "(3,4)"),
    (8*(a+b)^3, "(1,3)"),
    ((7/2)^3*(a+b+c)^3, "(1,4)"),
    ((5/2)^3*(b+c)^3, "(2,4)"),
    ((a+2*b)^3, "(1,2,3)"),
    ((2*a+b)^3, "(1,3,2)"),
    ((2*a+2*b+(7/2)*c)^3, "(1,3,4)"),
    (((7/2)*a+(7/2)*b+(3/2)*c)^3, "(1,4,3)"),
    ((b+(5/2)*c)^3, "(2,3,4)"),
    (((5/2)*b+(3/2)*c)^3, "(2,4,3)"),
    ((a+(7/2)*b+(7/2)*c)^3, "(1,2,4)"),
    (((7/2)*a+(5/2)*b+(5/2)*c)^3, "(1,4,2)"),
    ((a+(3/2)*c)^3, "(1,2)(3,4)"),
    ((2*a+(9/2)*b+(5/2)*c)^3, "(1,3)(2,4)"),
    (((7/2)*a+(9/2)*b+(7/2)*c)^3, "(1,4)(2,3)"),
    ((a+2*b+(7/2)*c)^3, "(1,2,3,4)"),
    ((a+(7/2)*b+(3/2)*c)^3, "(1,2,4,3)"),
    ((2*a+(9/2)*b+(7/2)*c)^3, "(1,3,2,4)"),
    ((2*a+b+(5/2)*c)^3, "(1,3,4,2)"),
    (((7/2)*a+(9/2)*b+(5/2)*c)^3, "(1,4,2,3)"),
    (((7/2)*a+(5/2)*b+(3/2)*c)^3, "(1,4,3,2)")
]

total_sum = 0
for coeff, cycle_str in terms:
    s = G(cycle_str)
    total_sum += coeff / theta_new(s)

result = total_sum.expand()

print("========================================")
print("After replacing s^{-1} with s, the simplified expression is:")
print(result)
print("========================================")
\end{verbatim}

\end{appendix}

\section*{Conflict of Interest Statement}
The authors declare that they have no known competing financial 
interests or personal relationships that could have appeared to 
influence the work reported in this paper.

\begin{Backmatter}

\end{Backmatter}

\end{document}